\documentclass{amsart}

\usepackage{amsthm,amsmath,amssymb,tikz,comment}

\title{Projections in the curve complex arising from covering maps}

\author{Robert Tang}
\address{Mathematics Institute, University of Warwick, Coventry, CV4 7AL, UK}
\email{robert.tang@warwick.ac.uk}
\urladdr{http://www.warwick.ac.uk/~mariam/}
\numberwithin{equation}{section}

\newtheorem{Thm}[equation]{Theorem}
\newtheorem{Prop}[equation]{Proposition}
\newtheorem{Lem}[equation]{Lemma}
\newtheorem{Cor}[equation]{Corollary}

\makeatletter
\newtheorem*{Thmmanual}{Theorem \@currentlabel}
\newenvironment{Thmman}[1]
  {\def\@currentlabel{#1}\Thmmanual}{\endThmmanual}
  
 \newtheorem*{Propmanual}{Proposition \@currentlabel}
\newenvironment{Propman}[1]
  {\def\@currentlabel{#1}\Propmanual}{\endPropmanual} 
\makeatother

\theoremstyle{definition}
\newtheorem{Dfn}[equation]{Definition}

\theoremstyle{remark}
\newtheorem{Rem}[equation]{Remark}
\newtheorem{Ex}[equation]{Example}

\newcommand{\R}{\mathbb{R}}
\newcommand{\Z}{\mathbb{Z}}

\newcommand{\Hy}{\mathbb{H}}

\newcommand{\TT}{\mathbb{T}}

\newcommand{\diam}{\mathrm{diam}}
\newcommand{\proj}{\mathrm{proj}}
\newcommand{\entry}{\mathrm{entry}}

\newcommand{\HD}{\mathrm{HausDist}}
\newcommand{\tvec}{\mathbf{t}}

\newcommand{\area}{\mathrm{area}}

\newcommand{\inte}{\mathrm{int}}
\newcommand{\rad}{\mathrm{rad}}
\newcommand{\Fix}{\mathrm{Fix}}
\newcommand{\cir}{\mathrm{circ}}

\newcommand{\norma}{\left\|\tvec\right\|_{\balpha}}
\newcommand{\norm}[1]{\left\|{#1}\right\|}
\newcommand{\short}{\mathrm{short}}
\newcommand{\Short}{\mathrm{Short}}

\newcommand{\width}{\mathrm{width}}
\newcommand{\length}{\mathrm{length}}

\newcommand{\Hull}{\mathrm{Hull}}
\newcommand{\Deck}{\mathrm{Deck}}
\newcommand{\Homeo}{\mathrm{Homeo}}
\newcommand{\one}{\mathbf{1}}

\newcommand{\Nb}{\mathcal{N}}
\newcommand{\C}{\mathcal{C}}
\newcommand{\X}{\mathcal{X}}
\newcommand{\Y}{\mathcal{Y}}
\newcommand{\A}{\mathcal{A}}

\newcommand{\Mod}{\mathrm{Mod}}

\newcommand{\bF}{\mathbf{F}}

\newcommand{\sK}{\mathsf{K}}
\newcommand{\sk}{\mathsf{k}}
\newcommand{\sQ}{\mathsf{Q}}

\newcommand{\sD}{\mathsf{D}}
\newcommand{\sW}{\mathsf{W}}
\newcommand{\sL}{\mathsf{L}}

\newcommand{\sr}{\mathsf{r}}
\newcommand{\sR}{\mathsf{R}}
\newcommand{\sLambda}{\mathsf{\Lambda}}

\newcommand{\rocket}{\mathfrak{R}}
\newcommand{\nose}{\mathfrak{N}}
\newcommand{\shaft}{\mathfrak{S}}
\newcommand{\fins}{\mathfrak{F}}

\newcommand{\balpha}{{\boldsymbol{\alpha}}}

\newcommand{\ang}[1]{\langle #1 \rangle}

\newcommand{\lab}[1]{\label{#1}}

\newcommand{\proofof}[1]{\hfill\newline\noindent\emph{Proof of {#1}.} }
\newcommand{\halmos}{\hfill$\square$}

\begin{document}

	\begin{abstract}
	Let $P : \Sigma \rightarrow S$ be a finite degree covering map between surfaces. Rafi and Schleimer show that there is an induced quasi-isometric embedding $\Pi : \C(S) \rightarrow \C(\Sigma)$ between the associated curve complexes. We define an operation on curves in $\C(\Sigma)$ using minimal intersection number conditions and prove that it approximates a nearest point projection to $\Pi(\C(S))$.
	We also approximate hulls of finite sets of vertices in the curve complex, together with their corresponding nearest point projections, using intersection numbers.
	\end{abstract}
	
		\maketitle
	
	\section{Introduction}\lab{chap-back}

\subsection{The curve complex}

	Let $S = (S,\Omega)$ denote a closed, orientable, connected surface of genus $g\geq 0$ together with a set $\Omega$ of $m\geq 0$ marked points. Define the \emph{complexity} of $S$ to be $\xi(S) := 3g - 3 + m$. 

	A \emph{curve} on $S$ is a continuous map $a : S^1 \rightarrow S - \Omega$, where $S^1 = \R/\Z$ is the circle. We will also write $a$ for its image on $S$. A curve $a$ is \emph{simple} if it is an embedded copy of $S^1$.

	We call a curve \emph{trivial} or \emph{peripheral} if it is freely homotopic to a curve bounding a disc or a disc with exactly one marked point respectively. A simple closed curve which is non-trivial and non-peripheral is called \emph{essential}. 	A \emph{multicurve} on $S$ is a finite collection of essential simple closed curves which can be realised disjointly simultaneously.

	Let $\C^0(S)$ denote the set of free homotopy classes of essential simple closed curves on $S$. Unless explicitly stated otherwise, we will blur the distinction between curves and their free homotopy classes.

	For what follows, we shall assume that $S$ has complexity $\xi(S)$ at least 2;  modifications to the following definition are required for low-complexity cases but we shall not deal with them here. For an introduction to the curve complex, see \cite{saul-notes}.
	
	\begin{Dfn}
	The \emph{curve complex} of $S$, denoted $\C(S)$, is a simplicial complex whose vertex set is $\C^0(S)$ and whose simplices are spanned by multicurves. In particular, two distinct simple closed curves are connected by an edge in $\C(S)$ if and only if they have disjoint representatives on $S$.
	\end{Dfn}
	
	The dimension of $\C(S)$ is equal to $\xi(S) - 1$. If one replaces marked points with boundary components then the top dimensional simplices of $\C(S)$ correspond to \emph{pants decompositions} of $S$, that is, multicurves which cut $S$ into a collection of pants.
	
	We endow $\C(S)$ with the standard simplicial metric: each $k$--simplex is isometrically identified with a standard Euclidean $k$--simplex whose edge lengths are equal to 1.	For our purposes, it suffices to study the 1--skeleton $\C^1(S)$ of the curve complex, often referred to as the \emph{curve graph}. Indeed $\C^1(S)$ equipped with the induced path metric, denoted $d_S$, is naturally \emph{quasi-isometric} to $\C(S)$. To simplify notation, we shall write $\C(S)$ in place of $\C^1(S)$ and $\alpha\in\C(S)$ to denote a curve (or multicurve).
	
	Given free homotopy classes of curves $\alpha$ and $\beta$, not necessarily simple, define their \emph{(geometric) intersection number} $i(\alpha,\beta)$ to be the minimal value of $|a \cap b|$ over all their representatives $a \in \alpha$ and $b \in \beta$ which are in general position on $S$.
	
	We say a finite collection of curves \emph{fills} $S$ if their complement is a disjoint union of discs each with at most one marked point.
	
	Let $\alpha$ and $\beta$ be curves in $\C(S)$. Their distance $d_S(\alpha,\beta)$ is equal to the length of a shortest edge-path in $\C(S)$ connecting $\alpha$ and $\beta$. Observe that $\alpha$ and $\beta$ are disjoint if and only if $d_S(\alpha,\beta) \leq 1$. We also have $d_S(\alpha,\beta) = 2$ if and only if $\alpha$ and $\beta$ intersect but do not fill $S$; and $d_S(\alpha,\beta) \geq 3$ if and only if they do fill.

	\begin{Lem}[\cite{hempel-cc}, \cite{saul-notes}]\lab{lemhempel}
	Suppose $\alpha$ and $\beta$ are curves in $\C(S)$. Then
	\[d_S(\alpha,\beta) \leq 2 \log_2 i(\alpha,\beta) + 2 \]
	whenever $i(\alpha,\beta) \neq 0$.\halmos
	\end{Lem}
	
	As an immediate corollary, we see that $\C(S)$ is connected (this was originally observed by Harvey in \cite{Harvey-cc}). It is also worth mentioning that one cannot give a lower bound on distance in $\C(S)$ purely in terms of intersection number -- indeed, one can find pairs of non-filling curves which intersect an arbitrarily large number of times.
	
	The curve graph is also locally infinite and has infinite diameter \cite{Kob-heights}. Masur and Minsky proved the following celebrated theorem regarding the large scale geometry of the curve graph:
	
	\begin{Thm}[\cite{MM1}]\lab{thmhypcc}
	Given any surface $S$ with $\xi(S) \geq 2$, there exists $\delta > 0$ so that the curve graph $\C(S)$ is $\delta$--hyperbolic. \halmos
	\end{Thm}
	
	In \cite{bhb-int}, Bowditch gives a combinatorial proof of hyperbolicity using intersection numbers. We will be extending many of the results established in his paper in Sections \ref{chap-hull} and \ref{chap-sing}.
	
	\begin{Thm}[\cite{bhb-unif}, \cite{Aougab-unif}, \cite{CRS-unif}, \cite{HPW-unicorn}]
	The constant $\delta > 0$ in Theorem \ref{thmhypcc} can be chosen independently of $S$. \halmos
	\end{Thm}
	
	Hensel, Przytycki and Webb in particular show that all geodesic triangles in $\C(S)$ possess 17--centres.
	
	\subsection{Statement of results}
	
		Let $P: \Sigma \rightarrow S$ be a finite degree covering map. For any simple closed curve $a\in\C(S)$, the preimage $P^{-1}(a)$ is a disjoint union of simple closed curves on $\Sigma$. Rafi and Schleimer define a (one-to-many) \emph{lifting operation} $\Pi : \C(S) \rightarrow \C(\Sigma)$ by setting $\Pi(a) = P^{-1}(a)$. Rafi and Schleimer first proved the following theorem using techniques from Teichm\"uller theory. In \cite{tang-covers}, we give a new proof using hyperbolic 3--manifold geometry.
		
		\begin{Thm}[\cite{saul-covers}]\lab{thmqiemb}
		 Let $P: \Sigma \rightarrow S$ be a finite degree covering map. Then the map $\Pi:\C(S) \rightarrow \C(\Sigma)$ defined above is a $\sLambda$--quasi-isometric embedding, where $\sLambda$ depends only on $\xi(\Sigma)$ and $\deg P$. \halmos
		\end{Thm}
		
		It follows that $\Pi(\C(S))$ is quasiconvex in $\C(\Sigma)$. Quasiconvexity is a particularly nice geometric property in $\delta$--hyperbolic spaces: for example, nearest point projections to quasiconvex subsets are coarsely well-defined (Lemma \ref{lemboundproj}).
		
		We define a combinatorial operation $\pi : \C(\Sigma) \rightarrow \Pi(\C(S)) \subseteq \C(\Sigma)$ as follows: Given a curve $\alpha\in\C(\Sigma)$, let $b$ be a curve which minimises $i(P(\alpha),\cdot)$ among all curves in $\C(S)$. We then set $\pi(\alpha) = \Pi(b)$. 
		The aim of this paper is to show that this operation satisfies the following coarse geometric properties. In particular, Theorem \ref{coverproj} can be viewed as a converse to Theorem \ref{thmqiemb}.
	
		\begin{Thmman}{\ref{coverproj}}
		 Let $P: \Sigma \rightarrow S$ be a finite degree covering map and suppose $\alpha\in\C(\Sigma)$ is a curve. Then $\pi(\alpha)$ is a uniformly bounded distance from any nearest point projection of $\alpha$ to $\Pi(\C(S))$ in $\C(\Sigma)$, where the bounds depend only on $\xi(\Sigma)$ and the degree of $P$.
		\end{Thmman}
		
		\begin{Propman}{\ref{covercirc}}
		 Assume further that $P$ is regular, with deck group $G$. Then $\pi(\alpha)$ is a uniformly bounded distance from any circumcentre for the $G$--orbit of $\alpha$ in $\C(\Sigma)$. Moreover, the bounds depend only on $\xi(\Sigma)$ and the degree of $P$.
		\end{Propman}

		In order to prove the above results, we develop combinatorial descriptions of hulls in the curve complex which may be of independent interest. We state simplified versions of the relevant propositions below -- see Section \ref{chap-hull} for more precise formulations.
		
		Let $\balpha = (\alpha_1, \ldots, \alpha_n)$ be an $n$--tuple of distinct curves in $\C(S)$, where $n\geq 2$. Given a non-zero vector $\tvec = (t_1, \ldots, t_n)$ of non-negative reals, let $\gamma_\tvec \in \C(S)$ be a curve which minimises the weighted intersection number $\sum_i t_i i(\alpha_i, \cdot)$. Define the \emph{hyperbolic hull} $\Hull(\balpha)$ to be the union of all geodesic segments in $\C(S)$ connecting a pair of points in $\balpha$ (viewed as a vertex set in $\C(S)$).
		
		\begin{Propman}{\ref{hullshort}}
		The sets $\Hull(\balpha)$ and $\bigcup_\tvec \gamma_\tvec$ agree up to a uniformly bounded Hausdorff distance in $\C(S)$, where the union is taken over all non-zero $\tvec \in \R_{\geq 0}^n$. Moreover, the bound depends only on $\xi(S)$ and $n$.
		\end{Propman}
		
		\begin{Propman}{\ref{proj}}
		Suppose $\beta\in\C(S)$ is a curve satisfying $i(\alpha_i,\beta) \neq 0$ for all $i$. Let $\tvec_\beta = (t_1,\ldots,t_n)$ be the vector given by $t_i = i(\alpha_i,\beta)^{-1}$ for each $i$. Then $\gamma_{\tvec_\beta}$ is a uniformly bounded distance from any nearest point projection of $\beta$ to $\Hull(\balpha)$ in $\C(S)$, where the bound depends only on $\xi(S)$ and $n$.
		\end{Propman}

	\subsection{Organisation}
	
	We begin by reviewing some coarse geometric notions in Section \ref{chap-coarse}, placing a particular emphasis on $\delta$--hyperbolic spaces. In Section \ref{chap-hull}, we introduce two notions of hulls for finite sets in $\C(S)$: one arising geometrically in $\C(S)$; the other defined using intersection number conditions. We give proofs of Propositions \ref{hullshort} and \ref{proj} assuming bounded diameter properties for sets of curves satisfying certain bounded weighted intersection numbers conditions (Lemma \ref{boundshort}) -- a key fact whose proof we defer to Section \ref{secproofshort}.
	In Section \ref{chap-cover2}, we utilise the results from Sections \ref{chap-coarse} and \ref{chap-hull} to give proofs of the main theorems. In Section \ref{chap-sing}, we generalise Bowditch's construction of singular Euclidean structures on surfaces \cite{bhb-int} on which the geodesic lengths of curves estimate suitable weighted intersection numbers. We show that these surface possess certain geometric properties, such as the existence of wide annuli, in order to control their sets of short curves and hence give a proof of Lemma \ref{boundshort}.

	\section*{Acknowledgements}

	I am grateful to both Brian Bowditch and Saul Schleimer for many helpful discussions and suggestions. I also thank Saul Schleimer and Marc Lackenby for thorough comments on my PhD thesis \cite{tang-thesis}, in which the contents of this paper are completely contained. This work was supported by a Warwick Postgraduate Research Scholarship.

	\section{Coarse geometry}\lab{chap-coarse}
	In this section, we recall some basic definitions and notions concerning Gromov hyperbolic spaces. Many of the statements and results are either well known in the literature or relatively straightforward to deduce; we shall include them for completeness and to establish notation and terminology. We refer the reader to \cite{BH-metric}, \cite{Gromov-hyp}, \cite{ABCFLMSS} and \cite{bhb-ggt} for more background.

\subsection{Notation}
	
		Let $(\X,d)$ be a metric space. Given any subset $A \subseteq \X$ and a point $x\in\X$, we define $d(x,A) := \inf \{d(x,a) ~|~ a\in A\}$. For $\sr \geq 0$, let
		\[\Nb_{\sr}(A) = \{x\in\X ~|~ d(x,A) \leq \sr \}\]
		denote the \emph{$\sr$--neighbourhood} of $A$ in $\X$.
		For subsets $A,B \subseteq \X$ and $\sr\geq 0$, write
		\[A \subseteq_\sr B \iff A \subseteq \Nb_\sr(B)\]
		 and
		\[A \approx_\sr B \iff A \subseteq_\sr B \textrm{ and } B \subseteq_\sr A.\]
		Define the \emph{Hausdorff distance} between $A$ and $B$ to be
		\[\HD(A,B) = \inf \{\sr \geq 0 ~|~ A \approx_\sr B\}. \]
			
		To simplify notation, we will often write $a\in\X$ in place of a singleton set $\{a\}\subseteq\X$.
				We will always use the standard Euclidean metric on the reals unless otherwise specified. If $a$ and $b$ are real numbers then
		\[a \approx_\sr b \iff |a-b| \leq \sr.\]
		We will also adopt the following notation:
						\[a \asymp_\sr b \iff a \leq \sr\, b + \sr ~\textrm{ and }~ b \leq \sr\, a + \sr.\]
		
		The \emph{diameter} of $A\subseteq\X$ is defined to be
		\[\diam(A) := \sup\{ d(x,y) ~|~ x,y\in A\}. \]

		\subsection{Geodesics, quasiconvexity and quasi-isometries}

		Let $I \subseteq \R$ be an interval. A \emph{geodesic} is a map $\gamma : I \rightarrow \X$ so that $d(\gamma(t),\gamma(s)) = |t-s|$ for all $t,s\in I$.
		A \emph{geodesic segment} connecting points $x$ and $y$ in $\X$ is the image of a geodesic $\gamma : [0,d(x,y)] \rightarrow \X$ such that $\gamma(0) = x$ and $\gamma(d(x,y)) = y$.
		A metric space $\X$ is called a \emph{geodesic space} if every pair of points can be connected by a geodesic segment. 
		
			A subset $U \subseteq \X$ is $\sQ$--\emph{quasiconvex} if any geodesic segment connecting any pair of points in $U$ lies in $\Nb_\sQ(U)$. We say a subset is \emph{quasiconvex} if it is $\sQ$--quasiconvex for some $\sQ\geq 0$.
			
			A (one-to-many) map $f: \X \rightarrow \Y$ between metrics spaces is a \emph{$\sLambda$--quasi-isometric embedding} if for all $x_1,x_2 \in \X$ and $y_1\in f(x_1)$, $y_2 \in f(x_2)$ we have
		\[d_{\Y}(y_1, y_2) \asymp_\sLambda d_{\X}(x_1,x_2).\]
		In addition, if $\Nb_\sLambda (f(\X)) = \Y$ then $f$ is called a \emph{$\sLambda$--quasi-isometry} and we say that $\X$ and $\Y$ are \emph{$\sLambda$--quasi-isometric}. If $\X$ and $\Y$ are $\sLambda$--quasi-isometric for some $\sLambda \geq 1$ then we may simply say that they are \emph{quasi-isometric}.

\subsection{Gromov hyperbolic spaces}

		We recall some basic results about Gromov hyperbolic spaces. Let $\X$ be a geodesic space.

		\begin{Dfn}[Geodesic triangle]
		A geodesic triangle $T$ in $\X$ consists of three points $x,y,z\in\X$ together with three geodesic segments $[x,y], [y,z], [z,x]$. The segments will be called the \emph{sides} of the triangle $T$.
		\end{Dfn}
		
		We will abbreviate $d(x,y)$ to $xy$.
				Let $x,y,z$ be points in $\X$. We will write
		\[ \ang{x,y}_z := \frac{1}{2}(xz + yz - xy). \]
		for the \emph{Gromov product} of $x$ and $y$ with respect to $z$.

		Given a geodesic triangle $T$ on $x,y,z\in\X$, we construct a \emph{comparison tripod} $\bar T$ as follows: Build a metric tree consisting of one central vertex whose valence is at most 3 and three vertices of valence one with three edges of lengths $\ang{y,z}_x, \ang{z,x}_y$ and $\ang{x,y}_z$. Label the central vertex $o_T$ and the other endpoints of the edges $\bar x, \bar y$ and $\bar z$ respectively. We allow for the possibility of edges having length zero in this construction.
		
		There exists a unique map
		\[\theta_T~:~T \rightarrow \bar T\]
		satisfying $\theta_T(x) = \bar x$, $\theta_T(y) = \bar y$ and $\theta_T(z) = \bar z$ which restricts to an isometric embedding on each edge of $T$. The elements of $\theta_T^{-1}(o_T)$ are called the \emph{internal points} of $T$.

		\begin{Dfn}[Thin triangle, $\delta$--hyperbolic space]
				A geodesic triangle $T$ is \emph{$\delta$--thin} if
		\[\diam (\theta_T^{-1}(p)) \leq \delta\]
		for all $p\in\bar T$. A geodesic space $\X$ is \emph{$\delta$--hyperbolic} if all of its geodesic triangles are $\delta$--thin. We call $\X$ \emph{(Gromov) hyperbolic} if it is $\delta$--hyperbolic for some $\delta \geq 0$.
		\end{Dfn}

		If $T$ is a $\delta$--thin geodesic triangle with vertices $x,y,z\in\X$ then its internal points decompose it into three pairs of $\delta$--fellow travelling geodesic segments whose lengths are $\ang{y,z}_x, \ang{z,x}_y$ and $\ang{x,y}_z$.

				The following result shows us that geodesic segments between two given points in a $\delta$--hyperbolic space are essentially unique up to bounded error.
		
		\begin{Lem}[Stability of geodesics]
		Let $x,y$ be points in a $\delta$--hyperbolic space $\X$. Then any two geodesic segments $\gamma_1, \gamma_2$ joining $x$ and $y$ $\delta$-fellow travel: if $u_1\in\gamma_1$ and $u_2\in\gamma_2$ are points such that $xu_1 = xu_2$ then $u_1u_2 \leq \delta$. In particular, geodesics in a $\delta$--hyperbolic space are $\delta$--quasiconvex		\end{Lem}
		
		\proof
		Consider the geodesic triangle with vertices $x$, $y$ and $y$ whose non-degenerate sides are $\gamma_1$ and $\gamma_2$. The result follows from the definition of $\delta$--thinness.
		\endproof

		We also state some equivalent notions of Gromov hyperbolicity:	
			
				\begin{Lem}[Four point condition, \cite{BH-metric} Proposition 1.22]\lab{4pt}
			If $\X$ is a $\delta$--hyperbolic space then
			\[xy + zw \leq \max \{ xz + yw,\, xw + yz \} + 2\delta\]
			for all $x,y,z,w\in\X$.
			
			Conversely, if the above inequality holds for all points $x,y,z$ and $w$ in a geodesic space $\X$, then $\X$ is $\delta'$--hyperbolic for some $\delta' \geq 0$ depending only on $\delta$. \halmos
			\end{Lem}
			
			 Suppose $\sk \geq 0$. A \emph{$\sk$--centre} for a geodesic triangle $T\subseteq\X$ is a point in $\X$ which lies within a distance $\sk$ of each side of $T$.

			\begin{Lem}[\cite{bhb-ggt} Proposition 6.13]
					Any geodesic triangle in a $\delta$--hyperbolic space possesses a $\delta$--centre, namely, any of its internal points.
			
			Conversely, suppose $\X$ is a geodesic space. If there is some $\sk\geq 0$ such that all geodesic triangles in $\X$ possess $\sk$--centres then $\X$ is $\delta$--hyperbolic for some $\delta\geq 0$ depending only on $\sk$. \halmos
			\end{Lem}

	\subsection{Nearest point projections to quasiconvex sets}

	Given a non-empty subset $U \subseteq \X$ and a point $x\in\X$, define
	\[\proj_U(x) := \{p \in U ~|~ xp = d(x,U) \}\]
	to be the set of \emph{nearest point projections} of $x$ to $U$ in $\X$. If $U$ is closed in $\X$ then $\proj_U(x)$ is always non-empty.

	Nearest point projections to geodesic segments can be approximated by internal points:

		\begin{Lem}\lab{lemoproj}
		Let $\X$ be a $\delta$--hyperbolic space. Let $T$ a geodesic triangle with vertices $x,y,z\in\X$. Let $o_x$ be the internal point of $T$ on $[y,z]$ and suppose $p\in [y,z]$ is a point such that $xp \leq d(x,[y,z]) + \epsilon$, for some $\epsilon \geq 0$. Then $o_xp \leq 2\delta + \epsilon$.
						\end{Lem}
		
		\proof
		Without loss of generality, suppose $p$ lies on $[o_x,y]$. Let $q$ be a point on $[x,y]$ such that $py = qy$. By $\delta$--hyperbolicity, we have $pq \leq \delta$. Let $o_z$ be the internal point on $[x,z]$ opposite $y$. Then
		\[xo_z + o_zq = xq \leq xp + \delta \leq d(x,[y,z]) + \epsilon + \delta \leq xo_x + \epsilon + \delta \leq xo_z + 2\,\delta + \epsilon\]
		and so $o_xp = o_zq \leq 2\delta + \epsilon$
		\endproof

	Let us now assume that $U$ is a closed, non-empty $\sQ$--quasiconvex subset of a $\delta$--hyperbolic space $\X$.
	
	\begin{Lem}\lab{lemop}
	 Let $p$ be a nearest point projection of $x\in\X$ to $U$. Let $u$ be any point in $U$ and let $o_x$, $o_p$ and $o_u$ be the respective internal points of a geodesic triangle with vertices $x$, $p$ and $u$. Then $po_x \leq \delta + \sQ$ and hence $po_u \leq \delta + \sQ$.
	\end{Lem}
	
	\proof
	By quasiconvexity of $U$ we have $d(o_x,U) \leq \sQ$. Thus,
	\[xo_u + o_up = xp = d(x,U) \leq xo_x + d(o_x,U) \leq xo_u + o_uo_x + \sQ \leq xo_u + \delta + \sQ\]
	and so $o_xp = o_up \leq \delta + \sQ$.
	\endproof
	
	\begin{Lem}\lab{lemboundproj}
	For all $x\in\X$,
	\[\diam(\proj_U(x)) \leq 2\delta + 2\sQ. \]
	\end{Lem}
	
	\proof
	Let $p$ and $q$ be nearest point projections of $x$ to $U$. Let $o_x$ be the respective internal point opposite $x$ of a geodesic triangle with vertices $x$, $p$ and $q$. Applying \ref{lemop}, we deduce $pq \leq po_x + o_xq \leq 2\delta + 2\sQ$.
	\endproof

A consequence of Lemma \ref{lemop} is that any geodesic from $x$ to a point in $U$ must pass within a distance of $\delta + \sQ$ of every nearest point projection of $x$ to $U$. It turns out that this property characterises nearest point projections to quasiconvex sets in hyperbolic spaces.
For $\sr \geq 0$, we define $\entry_U(x,\sr)$ to be the set of all points $q \in U$ such that for all $u\in U$, every geodesic connecting $x$ to $u$ passes within a distance of $\sr$ of $q$.
Such points will be called $\sr$--\emph{entry points} of $x$ to $U$.

	\begin{Lem}\lab{lemgeodentry2}	Let $\sr \geq 0$. Then for all $x\in\X$,
	\[\entry_U(x,\sr) \subseteq_{2\sr} \proj_U(x). \]	In particular, for $\sr \geq 2\delta + \sQ$ we have
	\[\entry_U(x,\sr) \approx_{2\sr} \proj_U(x).\]
	\end{Lem}
	
	\proof
		Suppose $p$ is a nearest point projection and $q$ is an $\sr$--entry point of $x$ to $U$ respectively. Then there is some point $y\in [x,p]$ so that $yq \leq \sr$. Now
	\[xy + yp = xp \leq xq \leq xy + yq \leq xy + \sr\]
	and so $pq \leq py + yq \leq 2\sr$ which proves the first statement. The second statement follows from Lemma \ref{lemop}.
	\endproof
	
	Furthermore, any geodesic segment from a point $x\in\X$ to $u\in U$ can be approximated by the concatenation of two segments: the first from $x$ to any point $p\in\proj_U(x)$ and the second from $p$ to $u$.
	
		\begin{Lem}\lab{lemgeodentry}
	Given $x\in\X$, let $p\in\proj_U(x)$. Then for any $u\in U$,
	\[ [x,u] \approx_{2\delta + \sQ} [x,p] \cup [p,u]\]
	and
	\[xu \approx_{2\delta + 2\sQ} xp + pu.\]
	\end{Lem}
	
	\proof
	By hyperbolicity, we have $[x,u] \approx_\delta [x,o_u] \cup [o_x,u]$ and $xu = xo_u + o_xu$. By Lemma \ref{lemop},
	\[\diam[o_u,p] = \diam[p,o_x] = po_x \leq \delta + \sQ\]
	and so $[o_u,p] \cup [p, o_x] \subseteq_{\delta + \sQ} \{o_u, o_x\} \subseteq_{\delta} [x,u]$. 	\endproof

	\subsection{Circumcentres}
			
		Let $\X$ be a $\delta$--hyperbolic space and suppose $U\subseteq \X$ is a non-empty finite subset.
		
				\begin{Dfn}[Radius, Circumcentre]
		The \emph{radius} of $U$ is
		\[\rad(U) := \min \left\{ \sr\geq 0 ~|~ \exists\, x\in\X, U \subseteq B_{\sr}(x) \right\}, \]
		where $B_\sr(x)$ is the closed ball of radius $\sr$ centred at $x$. We call a point $x\in\X$ a \emph{circumcentre} of $U$ if $U \subseteq B_{\sr}(x)$ for $\sr = \rad(U)$ and write $\cir(U)$ for the set of circumcentres of $U$.
		\end{Dfn}

		\begin{Lem}\lab{lemmiddist}
		Let $x$, $y$ and $z$ be points in $\X$. Suppose $m$ is a midpoint of some geodesic $[x,y]$ connecting $x$ and $y$. Then $\max \{xz, yz \} \approx_{\delta} \frac{1}{2}\, xy + mz$.
		\end{Lem}
		
		\proof
		Without loss of generality, suppose $xz \geq yz$. Then $m$ lies on $[x,o_z]$, where $o_z\in [x,y]$ is the internal point opposite $z$. By hyperbolicity, there is a point $q\in [x,z]$ such that $xq = xm$ and $qm\leq\delta$. Finally,
		\[\max \{xz, yz \} = xz = xq + qz \approx_\delta xm + mz = \frac{1}{2}\, xy + mz\]
		which completes the proof.
		\endproof
		
		\begin{Lem}\lab{lemcirc}
		Let $c$ be a circumcentre of $U$ and suppose $x\in\X$ is a point such that $U \subseteq B_{\sr + \epsilon}(x)$, where $\sr = \rad(U)$ and $\epsilon \geq 0$. Then $cx \leq 2\delta + 2\epsilon$ and hence $\diam(\cir(U))\leq 2\delta$.
		\end{Lem}
		
		\proof
		Let $m$ be a midpoint of $c$ and $x$. Choose $u \in U$ so that $um$ is maximal. Applying Lemma \ref{lemmiddist} and the definition of radius gives
		\[\rad(U) \leq um \leq \max\{cu, xu \} - \frac{1}{2}\, cx + \delta \leq \rad(U) +\epsilon - \frac{1}{2}\, cx + \delta\]
		and we are done.
		\endproof

		\begin{Lem}\lab{lemcirc2}
		Suppose $c$ is a circumcentre of $U$. Let $x$, $y\in U$ be points such that $xy \geq \diam(U) - 2\epsilon$, for some $\epsilon \geq 0$. Let $m$ be the midpoint of a geodesic segment $[x,y]$. Then $c \approx_{2\delta+\epsilon} m$. Furthermore, we have
		\[\diam(U) \leq 2\,\rad(U) \leq \diam(U) + 2\delta. \]
		\end{Lem}
		
		\proof
		Suppose $x'$ and $y'$ are points in $U$ satisfying $x'y' = \diam(U)$. By Lemma \ref{lemmiddist}, we deduce
		\[um \leq \max \{x'u, y'u\} - \frac{1}{2}\, x'y' + \delta \leq \diam(U) - \frac{1}{2}\,\diam(U) + \delta\]
		for all $u\in U$. Choosing $u$ so that $um$ is maximal yields
		\[\rad(U) \leq um \leq \frac{1}{2}\,\diam(U) + \delta.\]
		Next, observing
		\[\diam(U) = x'y' \leq x'c + cy' \leq 2\,\rad(U)\]
		completes the proof of the second claim. Finally, 
		\[cm \leq \max \{cx, cy\} - \frac{1}{2}\, xy + \delta \leq \rad(U) - \frac{1}{2}\,\diam(U) + \epsilon + \delta \leq 2\delta + \epsilon\]
		where we have applied the second claim and Lemma \ref{lemmiddist} once more.
		\endproof

	\section{Hulls in the curve complex}\lab{chap-hull}
	
Let $S = (S,\Omega)$ be a connected compact surface $S$ without boundary with a finite set of marked points $\Omega$ satisfying $\xi(S) \geq 2$. Throughout this section, we will fix an $n$-tuple ${\balpha = (\alpha_1, \ldots, \alpha_n)}$ of distinct multicurves in $\C(S)$, where $n\geq 2$. We will assume that no pair $\alpha_i$ and $\alpha_j$ has a common component.

We shall establish a coarse equality between two subsets of $\C(S)$ determined by $\balpha$ -- its hyperbolic hull $\Hull(\balpha)$, defined purely in terms of the geometry of $\C(S)$; and $\Short(\balpha,\sL)$ which is defined using only intersection numbers. We also give a combinatorial method of approximating nearest point projections to $\Hull(\balpha)$.

\subsection{Hyperbolic hulls}

	Let $\X$ be a $\delta$--hyperbolic space and suppose $U\subseteq \X$ is a set of points. The \emph{hyperbolic hull} of $U$, denoted $\Hull(U)$, is the union of all geodesic segments in $\X$ connecting a pair of points in $U$.

	\begin{Ex}
	 Let $U$ be a finite subset of $\Hy^n$, where $n\geq 1$. Then $\Hull(U)$ is a uniformly bounded Hausdorff distance away from the convex hull of $U$ in $\Hy^n$.
	\end{Ex}
	
	\begin{Lem}
	The hyperbolic hull of any non-empty set $U\subseteq \X$ is $2\delta$--quasiconvex.
	\end{Lem}
	
	\proof
	Let $u$ and $v$ be points in $\Hull(U)$. Let $x,y,z,w \in U$ be points, not necessarily distinct, so that $u \in [x,y]$ and $v \in [z,w]$. Let $[u,v]$ be any geodesic segment. By $\delta$--hyperbolicity, we have
	\[ [u,v] \subseteq_\delta [u,y] \cup [y,v] \subseteq_\delta [u,y] \cup [y,z] \cup [z,v] \subseteq [x,y] \cup [y,z] \cup [z,w] \subseteq \Hull(U),\]
	where $[u,y]$ and $[z,v]$ are assumed to be subarcs of $[x,y]$ and $[z,w]$ respectively.
	\endproof

	\begin{Lem}
	Suppose $C\subseteq \X$ is a $\sQ$--quasiconvex set which contains $U$. Then
	$\Hull(U) \subseteq_{\sQ} C$.
	
	\end{Lem}
	
	\proof
	This follows immediately from the definition of quasiconvexity.
	\endproof
	
	In fact, the above properties characterises hyperbolic hulls up to finite Hausdorff distance.
	
	\begin{Cor}
	Let $U\subseteq \X$ be non-empty. Suppose $C\subseteq \X$ is a $\sQ$--quasiconvex set such that
		\begin{enumerate}
			\item $C$ contains $U$, and
			\item for any $\sQ'$--quasiconvex set $C'\subseteq \X$ which also contains $U$, we have $C \subseteq_{\sr} C'$ for some $\sr = \sr(\sQ,\sQ') \geq 0$.
		\end{enumerate}
	Then $C$ and $\Hull(U)$ agree up to finite Hausdorff distance. \halmos
	\end{Cor}

\subsection{A hull via intersection numbers}
		
		\subsubsection{Short curve sets}

	Let $\balpha = (\alpha_1, \ldots, \alpha_n)$ be an $n$--tuple of distinct multicurves in $\C^0(S)$.
	A vector $\tvec = (t_1,\ldots,t_n) \neq \mathbf{0}$ of non-negative real numbers shall be referred to as a \emph{weight vector}. Write $\tvec\cdot\balpha$ for the formal sum $\sum_i t_i\alpha_i$. We will extend intersection number linearly over such sums:
	\[i(\tvec\cdot\balpha, \gamma) := \sum_i t_ii(\alpha_i,\gamma).\]
	For notational convenience, define a function on weight vectors by setting
	\[ \norma := \sqrt{i(\tvec\cdot\balpha)},\]
	where
	\[i(\tvec\cdot\balpha) = \sum_{j<k} t_jt_ki(\alpha_j,\alpha_k)\]
	is the \emph{self-intersection number} of $\tvec\cdot\balpha$. This serves as a rescaling factor for the singular Euclidean surface $S(\tvec\cdot\balpha)$ appearing in Section \ref{chap-sing}.

	Given $\sL \geq 0$, define 
		\[\short(\tvec\cdot\balpha, \sL) := \left\{\gamma\in\C(S) ~|~ i(\tvec\cdot\balpha,\gamma) \leq \sL \norma \right\}. \]
	If $\norma = 0$ then this set is contained in the 1--neighbourhood of $\balpha$.
	
	Note that $\short(\tvec\cdot\balpha, \sL)$ remains invariant under multiplying $\tvec$ by a positive scalar.
	
	When $\norma > 0$, the geodesic length of a curve $\gamma$ on $S(\tvec\cdot\balpha)$ approximates its intersection number with $\tvec\cdot\balpha$ (Proposition \ref{sing}). Thus, we can view $\short(\tvec\cdot\balpha, \sL)$ as the set of bounded length curves on $S(\tvec\cdot\balpha)$ rescaled to have unit area.

		\begin{Lem}\lab{boundshort}
		There exists a constant $\sL_0 > 0$ depending only on $\xi(S)$ such that, for any $\sL \geq \sL_0$, the set ${\short(\tvec\cdot\balpha, \sL)}$ is non-empty. Moreover,
		\[\diam_{\C(S)}(\short(\tvec\cdot\balpha, \sL)) \leq 4\log_2 \sL + \sk_0,\]
		where $\sk_0$ is a constant depending only on $\xi(S)$.
		\end{Lem}

		Consequently, up to bounded error, we can view $\short(\tvec\cdot\balpha, \sL)$ as a single curve in $\C(S)$ which has minimal intersection number with $\tvec\cdot\balpha$. The proof of the above lemma will be given in Section \ref{secproofshort} and largely follows the arguement in Bowditch's paper (\cite{bhb-int} Lemma 4.1).

		\subsubsection{The short curve hull}

		For $\sL\geq 0$, define the $\sL$\emph{--short curve hull} of $\balpha$ to be
		\[\Short(\balpha,\sL) := \bigcup_{\tvec} ~\short(\tvec\cdot\balpha, \sL),\]
		where the union is taken over all weight vectors $\tvec \in\R_{\geq 0}^{n}$ (or, equivalently, by choosing one representative from each projective class).
		
		We write $\Hull(\balpha)\subseteq\C(S)$ for the hyperbolic hull of $\balpha$ considered as a set of vertices in $\C(S)$.
		
		\begin{Prop}\lab{hullshort}
		Let $\balpha$ be an $n$--tuple of multicurves in $\C(S)$. Then for any $\sL\geq \sL_0$,
		\[\Short(\balpha, \sL) \approx_{\sk_1} \Hull(\balpha)\]
		where $\sk_1$ depends only on $\xi(S)$, $n$ and $\sL$.
		\end{Prop}

		This is essentially an extension of Bowditch's coarse description of geodesics using intersection numbers employed in his proof of hyperbolicity of the curve complex \cite{bhb-int}. Let us begin with a reformulation of his result:
		
		\begin{Lem}[\cite{bhb-int} Proposition 6.2]\lab{bowdlines}
		Let $\balpha' = (\alpha_1, \alpha_2)$ be a pair of multicurves in $\C(S)$. Let $[\alpha_1,\alpha_2]$ denote any geodesic segment connecting $\alpha_1$ and $\alpha_2$ in $\C(S)$. Then for all $\sL \geq \sL_0$, we have
		\[ \Short(\balpha', \sL) \approx_{\sk'_1} [\alpha_1,\alpha_2].\]
		where $\sk'_1\geq 0$ depends only $\xi(S)$ and $\sL$. \halmos
		\end{Lem}

		\proofof{Proposition \ref{hullshort}}
		By applying the previous lemma to all pairs of multicurves $(\alpha_i,\alpha_j)$ in $\balpha = (\alpha_1,\ldots,\alpha_2)$, we obtain the inclusion:
		\[\Hull(\balpha) \subseteq_{\sk'_1}\Short(\balpha, \sL). \]

		Fix a weight vector $\tvec=(t_1,\ldots,t_n)$ and assume, for notational simplicity, that the quantity $t_jt_ki(\alpha_j,\alpha_k)$ is maximised when $\{j,k\} = \{1,2\}$. Let ${\balpha' = (\alpha_1,\alpha_2)}$ and $\tvec' = (t_1,t_2)$. Since there are $\frac{n(n-1)}{2}$ distinct unordered pairs of indices $\{j,k\}$, it follows that
		\[\norma^2 = \sum_{j<k} t_jt_ki(\alpha_j,\alpha_k) \leq \frac{n(n-1)}{2} t_1t_2i(\alpha_1,\alpha_2) = \frac{n(n-1)}{2} \norm{\tvec'}_{\balpha'}^2.\]

	Let $\gamma$ be a curve in $\short(\tvec\cdot\balpha, \sL)$. 

	Then
		\[i(\tvec'\cdot\balpha',\gamma) \leq i(\tvec\cdot\balpha,\gamma) \leq \sL\norma \leq \sL\sqrt{\frac{n(n-1)}{2}\norm{\tvec'}_{\balpha'}^2} \leq \frac{n\sL}{\sqrt{2}}\norm{\tvec'}_{\balpha'}\]

		which implies
		\[\short(\tvec\cdot\balpha, \sL) \subseteq \short\left(\tvec'\cdot\balpha', \frac{n\sL}{\sqrt{2}}\right).\]

		Invoking Lemma \ref{bowdlines}, we have
		\[ \short\left(\tvec'\cdot\balpha', \frac{n\sL}{\sqrt{2}}\right) \subseteq_{\sr} [\alpha_1,\alpha_2] \subseteq \Hull(\balpha)\]
		where $\sr\geq 0$ is some constant depending on $n$, $\sL$ and $\xi(S)$.

		This concludes the proof of Proposition \ref{hullshort}.
		\halmos

		We can describe the above proof in terms of the geometry of $S(\tvec\cdot\balpha)$ (see Section \ref{chap-sing} for the construction). Assume $S(\tvec\cdot\balpha)$ has unit area. One can obtain $S(\tvec'\cdot\balpha')$ by homotoping the annuli consisting of rectangles traversed by $\alpha_i$ to the core curve $\alpha_i$ for each $i \neq 1,2$. The maximality assumption on $\alpha_1$ and $\alpha_2$ ensures that the total area of the remaining rectangles is at least $\frac{2}{n(n-1)}$. We then scale $S(\tvec'\cdot\balpha')$ by a factor of at most $\frac{n}{\sqrt{2}}$ to give it unit area. Finally, observe that the length of a curve $\gamma$ on $S(\tvec\cdot\balpha)$ can only increase by a factor of at most $\frac{n}{\sqrt{2}}$ during this process.
		
		\subsection{Nearest point projections to hulls}
		
		In this section, we approximate nearest point projections to short curve hulls using only intersection number conditions.

		\begin{Dfn}
		Let $\beta\in\C(S)$ be a multicurve. A weight vector $\tvec = (t_1,\ldots,t_n)$ satisfying
			\[t_ji(\alpha_j,\beta) = t_ki(\alpha_k,\beta)\]
			for all $j,k$ is called a \emph{balance vector} for $\beta$ with respect to $\balpha$.
		\end{Dfn}	
			If $\beta$ intersects all $\alpha_i$ then setting $t_i = i(\alpha_i,\beta)^{-1}$ yields the unique balance vector up to positive scale. If not, we can set $t_i = 1$ whenever $i(\alpha_i,\beta) = 0$ and $t_i = 0$ otherwise to produce a balance vector. We will use $\tvec_\beta$ to denote any balance vector for $\beta$. We also remark that the above definition is analogous to the notion of \emph{balance time} for quadratic differentials as described by Masur and Minsky \cite{MM1}.

		\begin{Prop}\lab{proj}
		Assume $\sL \geq \sL_0$. 
		Given a multicurve $\beta\in\C(S)$, let $\tvec_\beta$ be any balance vector with respect to $\balpha$. Let $\gamma$ be any nearest point projection of $\beta$ to $\Hull(\balpha)$. Then
		\[\gamma \approx_{\sk_2} \short(\tvec_\beta\cdot\balpha, \sL),\]
		where $\sk_2\geq 0$ depends only on $\xi(S)$, $n$ and $\sL$.

		\end{Prop}

		As was the case with Proposition \ref{hullshort}, this is an extension of a result of Bowditch. His result was originally phrased in terms of centres for geodesic triangles, however, our statement agrees with it up to uniformly bounded error.
		
		\begin{Lem}[\cite{bhb-int} Proposition 3.1 and Section 4]\lab{centre}
		Let $\alpha_1$, $\alpha_2$ and $\beta$ be multicurves in $\C(S)$. Let $\tvec'_\beta$ be a balance vector for $\beta$ with respect to $\balpha' = (\alpha_1,\alpha_2)$. Let $\gamma$ be a nearest point projection of $\beta$ to $[\alpha_1,\alpha_2]$. Then
		\[\gamma \approx_{\sk'_2} \short(\tvec'_\beta\cdot\balpha', \sL),\]
		where $\sk'_2$ depends only on $\xi(S)$ and $\sL$. \halmos
		\end{Lem}
		
		If $\beta$ is disjoint from some $\alpha_i$ then Proposition \ref{proj} follows immediately from Hempel's bound (Lemma \ref{lemhempel}). We will henceforth assume this is not the case.
		
		Our first step is to reduce the problem of finding a nearest point projection to a hyperbolic hull to that of projecting to a suitable geodesic. 

		\begin{Lem} \lab{projhull}
		Let $U$ be a subset of a $\delta$--hyperbolic space $\X$. Fix a point $w\in\X$. Assume there exist $x,y\in U$ and $\sR\geq 0$ such that
		\[d_{\X}([x,y],[z,w]) \leq \sR.\]
		for all $z\in U$. Let $p$ and $q$ be nearest point projections of $w$ to $\Hull(U)$ and $[x,y]$ respectively. Then
		\[p \approx_{\sR'} q \]
		where $\sR'$ depends only on $\sR$ and $\delta$.
		\end{Lem}

		\proof
		By Lemma \ref{lemgeodentry2}, it suffices to show that for all $u\in\Hull(U)$, any geodesic $[w,u]$ must pass within a bounded distance of $q$. If $u$ lies on a geodesic segment $[z,z']$ for some $z,z'\in U$ then $[w,u]$ must lie inside the $2\delta$--neighbourhood of $[w,z]$ or $[w,z']$. Hence, we only need to bound $d(q, [w,z])$ for all $z\in U$ in terms of $\delta$ and $\sR$.
		
		Recall that geodesic segments are $\delta$--quasiconvex. If $z$ coincides with $x$ or $y$ then $d(q, [w,z]) \leq 3\delta$ by Lemma \ref{lemop}. Let us now suppose that $x$, $y$ and $z$ are distinct. Choose points $u\in [z,w]$ and $v\in [x,y]$ so that $uv = d_{\X}([x,y],[z,w]) \leq \sR$. Then
		\[q \subseteq_{3\delta} [w,v] \subseteq_{\sR + \delta} [w,u] \subseteq [w,z]\]
		where we have applied Lemma \ref{lemgeodentry} for the first comparison.

		\endproof

		In order to exploit the above result, we recall yet another lemma of Bowditch:

		\begin{Lem}[\cite{bhb-int} Proposition 6.3]\lab{4ptint}
		Suppose ${\alpha_1,\alpha_2,\alpha_3,\alpha_4\in\C(S)}$ are multicurves which satisfy
		\[i(\alpha_1,\alpha_4)\:i(\alpha_2,\alpha_3) \leq \sr~i(\alpha_1,\alpha_2)\;i(\alpha_3,\alpha_4)\]
		for some $\sr>0$. Then
		\[ d_S([\alpha_1,\alpha_2] , [\alpha_3,\alpha_4]) \leq \sR\]
		where $\sR\geq 0$ depends only on $\sr$ and $\xi(S)$. \halmos
		\end{Lem}

		\proof[Proof of Proposition \ref{proj}]
		Let $\tvec_\beta$ be a balance vector for $\beta$ with respect to $\balpha$. To simplify notation, assume $t_jt_ki(\alpha_j,\alpha_k)$ is maximised when $\{j,k\} = \{1,2\}$. Let $\gamma_{12}$ and $\gamma_0$ be nearest point projections of $\beta$ to $[\alpha_1,\alpha_2]$ and $\Hull(\balpha)$ respectively.

		This implies
		\[t_2t_ji(\alpha_2,\alpha_j) \leq t_1t_2i(\alpha_1,\alpha_2)\]
		for any $j = 1, \ldots, n$. As $\beta$ is assumed to intersect all the $\alpha_i$, we have ${t_i = i(\alpha_i,\beta)^{-1}}$ (after rescaling) and so
		\[i(\alpha_1,\beta)\;i(\alpha_2,\alpha_j) \leq i(\alpha_1,\alpha_2)\;i(\alpha_j,\beta).\]
		Invoking Lemma \ref{4ptint} gives
		\[ d_S([\alpha_1,\alpha_2] , [\alpha_j,\beta]) \leq \sR\]
		and so by Lemma \ref{projhull} we deduce
		\[d_S(\gamma_{12},\gamma_0) \leq \sR',\]
		where $\sR'$ depends only on $\xi(S)$.

	Now suppose $\gamma$ is a curve in $\short(\tvec_\beta\cdot\balpha, \sL)$. Using the same reasoning as for the proof of Proposition \ref{hullshort}, we see that
	\[\gamma\in\short(\tvec_\beta\cdot\balpha, \sL) \subseteq \short\left(\tvec'_\beta\cdot\balpha', \frac{n\sL}{\sqrt{2}}\right)\]
	where $\balpha' = (\alpha_1,\alpha_2)$ and $\tvec'_\beta = (t_1,t_2)$. By Lemma \ref{centre}, we deduce that
	\[d_S(\gamma,\gamma_{12}) \leq \sk'_2,\]
	for some $\sk'_2$ depending only on $n$, $\sL$ and $\xi(S)$. This together with the preceding inequality implies
	\[d_S(\gamma,\gamma_0) \leq \sR' + \sk'_2\]
	which concludes the proof of the proposition.
	\endproof		
	
	\section{Covering maps}\lab{chap-cover2}

	\subsection{Operations on curves arising from covering maps}\lab{secprojint}
	
		We first recall some definitions and notation. 
		Let $P: \Sigma \rightarrow S$ be a finite degree covering map of surfaces. The preimage $P^{-1}(a)$ of a simple closed curve $a$ on $S$ under $P$ is a multicurve on $\Sigma$. This induces a one-to-many \emph{lifting map} ${\Pi : \C(S) \rightarrow \C(\Sigma)}$ between curve complexes by setting $\Pi(a) := P^{-1}(a)$. Recall the following theorem of Rafi and Schleimer:

		\begin{Thmman}{\ref{thmqiemb}}[\cite{saul-covers}]
		 Let $P: \Sigma \rightarrow S$ be a finite degree covering map. Then the map $\Pi:\C(S) \rightarrow \C(\Sigma)$ defined above is a $\sLambda$--quasi-isometric embedding, where $\sLambda$ depends only on $\xi(\Sigma)$ and $\deg P$.
		\end{Thmman}

		It immediately follows that $\Pi(\C(S))$ is quasiconvex in $\C(\Sigma)$. This naturally leads to the question of understanding nearest point projections to $\Pi(\C(S))$. Let us define an operation $\pi: \C(\Sigma) \rightarrow \Pi(\C(S))$ as follows: given a curve $\alpha\in\C(\Sigma)$, let $b\in\C(S)$ be a curve which has minimal intersection number with $P(\alpha)$ on $S$ and set $\pi(\alpha) = \Pi(b)$. We will prove the following in Section \ref{secproj}

		\begin{Thm}\lab{coverproj}
		Let $P: \Sigma \rightarrow S$ be a finite degree covering map of surfaces and let $\Pi$ and $\pi$ be as above. Given a curve $\alpha\in\C(\Sigma)$, let $\gamma$ be a nearest point projection of $\alpha$ to $\Pi(\C(S))$ in $\C(\Sigma)$. Then $\pi(\alpha) \approx_{\sk_3} \gamma$, where $\sk_3$ is a constant depending only on $\deg P$ and $\xi(\Sigma)$.
		\end{Thm}

		Consequently, the operation $\alpha \mapsto \pi(\alpha)$ is coarsely well-defined.
		One can check that the minimal value of $i(P(\alpha),\cdot)$ over all  closed curves on $S$ is attained by a simple closed curve. In Section \ref{seccirc}, we will also prove the following.

		\begin{Prop}\lab{covercirc}
		Suppose further that $P$ is regular and let $G$ be its group of deck transformations. Let $\gamma'$ be a circumcentre of the $G$--orbit of a curve $\alpha$ in $\C(\Sigma)$. Then $\pi(\alpha) \approx_{\sk_4} \gamma'$, where $\sk_4$ is some constant depending only on $\deg P$ and $\xi(\Sigma)$.
		\end{Prop}

		Recall that the \emph{deck transformation group} $\Deck(P)$ of a covering map $P: \Sigma \rightarrow S$ is the group of all homeomorphisms $f \in \Homeo(\Sigma)$ satisfying $P \circ f = P$. In order for the above statement to make sense, we must check that $\Deck(P)$ can be identified with its image in the \emph{mapping class group} $\Mod(\Sigma) = \Homeo(\Sigma)/\Homeo_0(\Sigma)$.
		
		\begin{Lem}
		Let $P: \Sigma \rightarrow S$ be a finite degree covering map between surfaces of negative Euler characteristic. Then the natural quotient map from $\Deck(P)$ to $\Mod(\Sigma)$ is injective.
		\end{Lem}

		\proof
		We will only give a sketch proof. Endow $\inte (S)$ with a hyperbolic metric and pull it back to $\inte (\Sigma)$ via $P$. The group $\Deck(P)$ then acts on $\inte(\Sigma)$ by isometries. The result follows since any isometry of a hyperbolic surface isotopic to the identity must in fact coincide with the identity. \endproof

		Note, however, that the above lemma does not hold for covers of the torus or annulus.

		\subsection{Nearest point projections}\lab{secproj}
	
		\subsubsection{Regular covers}
		
		We shall first deal with the case where $P: \Sigma \rightarrow S$ is a regular cover. Let $G \leq \Mod(\Sigma)$ be the group of deck transformations of $P$. Given a curve $\alpha\in\C(\Sigma)$, observe that the set of lifts of $P(\alpha)$ to $\Sigma$ via $P$ is exactly $G\alpha$. Let $\balpha = (\alpha_1, \ldots, \alpha_n)$ be an $n$--tuple of curves whose entries are the lifts of $P(\alpha)$ in any order. Note that $n\geq 1$ is some divisor of $\deg P$. Let $\one$ denote the vector of length $n$ with all entries equal to 1.

		\begin{Lem}\lab{lempishort}
		Let $\alpha$ and $\balpha$ be as above. Then $\pi(\alpha) \in \short(\one\cdot\balpha, \sL_0 |G|)$ where $\sL_0$ is a constant depending only on $\xi(\Sigma)$.
		\end{Lem}

		\proof
		Let $b$ be a closed curve on $S$. Each point of ${b \cap P(\alpha)}$ on $S$ lifts to exactly $|G| = \deg P$ points of ${P^{-1}(b) \cap G\alpha}$ on $\Sigma$ via $P$, hence
		\[i(P^{-1}(b), \balpha) = |G|~ i(b, P(\alpha)).\]
		By Lemma \ref{boundshort}, there exists a curve $\gamma\in\C(\Sigma)$ such that
		\[i(\gamma, \alpha) \leq \sL_0 \norm{\one}_{\balpha}\]
		for some constant $\sL_0 = \sL_0(\xi(\Sigma))$.
		Now assume $b$ has minimal intersection with $P(\alpha)$ out of all curves on $S$. It follows that
		\[i(P^{-1}(b), \balpha) = |G|~ i(b, P(\alpha)) \leq |G|~ i(P(\gamma), P(\alpha)) = i(G\gamma,\balpha) \leq |G|~ i(\gamma, \balpha).\]
		Finally, by combining the preceding inequalities, we see that
		\[i(\pi(\alpha),\balpha) = i(P^{-1}(b), \balpha) \leq |G|~ i(\gamma, \balpha) \leq |G|~ \sL_0 \norm{\one}_{\balpha}.\]
		Thus $\pi(\alpha) \in \short(\one\cdot\balpha, \sL)$ for $\sL = \sL_0 |G|$.
		\endproof

		\begin{Lem}\lab{lemprojmincurve}
		Let $\gamma$ be any curve in $\Pi(\C(S))$ and let $\beta$ be any of its nearest point projections to $\Hull(\balpha)$. Then $d_{\Sigma}(\pi(\alpha),\beta) \leq \sk_4$, where $\sk_4$ depends only on $\deg P$ and $\xi(\Sigma)$.
		\end{Lem}

		\proof
		We may replace $\gamma$ with the multicurve $G\gamma$ since their nearest point projections to $\Hull(\balpha)$ are a uniformly bounded distance apart. Since $G$ acts transitively on $G\alpha$ and leaves $G\gamma$ invariant, it follows that $i(G\gamma,\alpha_i) = i(G\gamma,\alpha_j)$ for all $i,j$. Thus, $\one$ serves as a balance vector for $G\gamma$ with respect to $\balpha$. By Proposition \ref{proj}, we deduce that
		\[\beta \approx_{\sk_2} \short(\one\cdot\balpha, \sL),\]
		where $\sk_2$ depends only on $\xi(\Sigma)$, $n$ and $\sL\geq \sL_0$. Applying the previous lemma completes the proof.
		\endproof

		\proofof{Theorem \ref{coverproj} for regular covers}
		Let $\alpha$ and $\balpha$ be as above. Let $\gamma$ be any curve in $\Pi(\C(S))$. Since $\Hull(\balpha)$ is quasiconvex, Lemmas \ref{lemprojmincurve} and \ref{lemgeodentry} imply that any geodesic connecting $\alpha$ to $\gamma$ in $\C(\Sigma)$ must pass within a distance $\sr$ of $\pi(\alpha)$, where $\sr$ depends only on $\deg P$ and $\xi(\Sigma)$. Therefore $\pi(\alpha)$ is an $\sr$--entry point of $\alpha$ to $\Pi(\C(S))$. Since $\Pi(\C(S))$ is also quasiconvex, Lemma \ref{lemgeodentry2} implies $\pi(\alpha)$ is a uniformly bounded distance away from any nearest point projection of $\alpha$ to $\Pi(\C(S))$.

		\halmos
		
		\subsubsection{The general case}

		The main obstacle in proving Theorem \ref{coverproj} for a non-regular covering map $P: \Sigma \rightarrow S$ is the following: given a simple closed curve $\alpha\in\C(\Sigma)$ there may be some lifts of $P(\alpha)$ to $\Sigma$ which are not simple. To address this issue, we pass to a suitable finite cover of $\Sigma$ using a standard group theoretic argument.

		\begin{Lem}
		Let $P: \Sigma \rightarrow S$ be a covering map of finite degree. Then there exists a cover $Q : \hat{\Sigma} \rightarrow \Sigma$ such that $F := P \circ Q$ is regular and $\deg F \leq (\deg P)!$. 
		\end{Lem}

		\proof
		Let $H$ be the finite index subgroup of $\Gamma = \pi_1(S)$ corresponding to the covering map $P$ and let $H_0$ be the intersection of all $\Gamma$--conjugates of $H$. It is straightforward to check that $H_0$ is exactly the kernel of the action of $\Gamma$ on the set of left cosets of $H$ by left-multiplication. The desired result then follows.
		\endproof

		The covering map $F$ defined above is universal in the sense that any regular cover of $S$ which factors through $P$ must also factor through $F$.

		\begin{Lem}
		Let $P: \Sigma \rightarrow S$ and $F : \hat{\Sigma} \rightarrow S$ be as above. If $\alpha$ is a simple closed curve on $S$ then all lifts of $P(\alpha)$ to $\hat \Sigma$ via $F$ are simple.
		\end{Lem}

		\proof
		Any lift of $\alpha$ to $\hat \Sigma$ via $Q$ is also a simple lift of $P(\alpha)$ via $F$. Since $F$ is regular, it follows that all other lifts of $P(\alpha)$ to $\hat \Sigma$ are simple.
		\endproof

		Before continuing with the proof, we first show that that nearest point projections to quasiconvex sets are well-behaved under quasi-isometric embeddings. We remind the reader that we allow $f$ to be a one-to-many function.
		
		\begin{Lem}
		Let $f:\X \rightarrow \X'$ be a $\sLambda$--quasi-isometric embedding of geodesic spaces, where $\X'$ is $\delta'$--hyperbolic. Let $C$ be a $\sQ$--quasiconvex subset of $\X$ and let $C' = f(C)$. Given a point $x\in\X$, let $x'$ be a point in $f(x)$. Let $p$ and $q'$ be nearest point projections of $x$ to $C$ and $x'$ to $C'$ respectively. Let $q\in\X$ be a point so that $q'\in f(q)$. Then $p \approx_{\sK} q$, where $\sK$ depends only on $\delta'$, $\sLambda$ and $\sQ$.
		\end{Lem}

		\proof
		First, note that $\X$ is $\delta$--hyperbolic and $C'$ is $\sQ'$--quasiconvex in $\X'$ for some constants $\delta = \delta(\sLambda,\delta')$ and $\sQ' = \sQ'(\sQ,\sLambda,\delta)$. Let $c\in\X$ be a $\sk$--centre for $x$, $p$ and $q$, where $\sk=\delta$. Any point $c' \in f(c)$ is then a $\sk'$--centre for $x'$, $p'$ and $q'$, where $\sk' = \sk'(\sk,\sLambda)$ and $p'\in f(p)$. One can check that $xp \approx_{2\sk} xc + cp$. By quasiconvexity of $C$, there is some point $y\in C$ satisfying $cy \leq \sk + \sQ$. Since $p$ is a nearest point projection of $x$ to $C$, we obtain
		\[xc + cp -2\sk \leq xp \leq xy \leq xc + cy \leq xc + \sk + \sQ\]
		which implies $cp \leq \sQ + 3\sk$. Similarly, we can deduce $c' q' \leq \sQ' + 3 \sk'$. Since $f$ is a $\sLambda$--quasi-isometric embedding, it follows that $cq \leq \sLambda\times c' q' + \sLambda$ and hence
		\[pq \leq pc + cq \leq \sK,\]
		where $\sK = \sQ + 3\sk + \sLambda(\sQ' + 3 \sk') + \sLambda$.
		\endproof

		Let $\Phi : \C(S) \rightarrow \C(\hat \Sigma)$ and $\Psi : \C(\Sigma) \rightarrow \C(\hat \Sigma)$ be the lifting maps induced by the covering maps $F$ and $Q$ respectively. Define $\phi : \C(\hat \Sigma) \rightarrow \Phi(\C(S))$
		
		to be the projection map associated to $F$ as described in Section \ref{secprojint}. We may assume that $\phi \circ \Psi= \Psi \circ \pi$. 

		\proofof{Theorem \ref{coverproj}}
		Given $\alpha\in\C(\Sigma)$, let $\hat \alpha$ be any of its lifts to $\hat\Sigma$ via $Q$. Note that $\phi(\hat\alpha) = \Psi(\pi(\alpha))$. Let $\hat\gamma$ be a nearest point projection of $\hat\alpha$ to $\Phi(\C(S))$ in $\C(\hat\Sigma)$ and let $\gamma = Q(\hat\gamma) \in \Pi(\C(S))$.
		Since $F$ is regular, we can apply Theorem \ref{coverproj} to deduce that
		\[d_{\hat\Sigma}(\phi(\hat\alpha),\hat\gamma) \leq \hat \sk_3,\]
		where $\hat \sk_3$ depends only on $\deg F$ and $\xi(\hat\Sigma)$ which can in turn be bounded in terms of $\deg P$ and $\xi(\Sigma)$. By Theorem \ref{thmqiemb}, $\Psi$ is a $\sLambda$--quasi-isometric embedding, where $\sLambda = \sLambda(\deg F,\xi(\hat\Sigma))$, and so 
		\[d_{\Sigma}(\pi(\alpha),\gamma) \leq \sLambda\hat \sk_3 + \sLambda.\]
		By the previous lemma, $\gamma$ is a uniformly bounded distance away from any nearest point projection of $\alpha$ to $\Pi(\C(S))$ in $\C(\Sigma)$ and we are done.
		\halmos

		\subsection{Circumcentres of orbits}\lab{seccirc}

		We now show that $\pi(\alpha)$ also approximates a circumcentre of $G \alpha$ in $\C(\Sigma)$, where $G$ is the deck transformation group of a regular cover $P : \Sigma \rightarrow S$. First, we give the following characterisation of circumcentres of orbits under finite group actions on $\delta$--hyperbolic spaces:

		\begin{Lem}
		Assume $G$ is a finite group acting by isometries on a $\delta$--hyperbolic space $\X$. Fix a point $x_0\in\X$ and let $c$ be a circumcentre for $Gx_0$. Given a point $z\in\X$, let $p$ be any of its nearest point projection to $\Hull(Gx_0)$.

		Then
		\[pc \leq \rad(Gz) + 7\delta\]
		and hence
		 \[zc \leq \rad(Gz) + d(z,\Hull(Gx_0)) + 7\delta. \]
		\end{Lem}
		
		\proof
		We first claim that $p$ lies within a distance $\delta$ of a geodesic segment $[u,v]$, where $u,v\in Gx_0$ are points such that $uv \geq \diam(Gx_0) - 2\delta$. Suppose $p$ lies on a geodesic segment $[x,y]$ for some $x$ and $y$ in $Gx_0$. There exist some $x'$ and $y'$ in $Gx_0$ so that $xx' = yy' = \diam(Gx_0)$. If $x' = y '$ then the claim follows from hyperbolicity. Now assume $x' \neq y'$. By Lemma \ref{4pt}, we have
		\[2\,\diam(Gx_0) = xx' + yy' \geq \max \{xy + x'y', xy' + x'y \} \geq 2\,\diam(Gx_0) - 2\delta.\]
		If $xy + x'y' \geq 2\,\diam(Gx_0) - 2\delta$ then $xy \geq \diam(Gx_0) - 2\delta$ which implies the claim. If not then $xy' \geq \diam(Gx_0) - 2\delta$. The claim then follows by considering a geodesic triangle with $x$, $y$ and $y'$ as its vertices.

		Now suppose $q\in[u,v]$ is a point so that $pq \leq \delta$. Then
		\[d(z,[u,v]) \leq zq \leq zp + pq \leq d(z,\Hull(Gx_0)) + \delta \leq d(z,[u,v]) + \delta. \]
		Invoking Lemma \ref{lemoproj} we obtain $p \approx_\delta q \approx_{3\delta} o$, where $o\in[u,v]$ is the internal point opposite $z$. Setting $\sD := \rad(Gz) \geq \frac{1}{2}\,\diam(Gz)$, observe that
		\[d(z,x_0) = d(gz,gx_0) \approx_{2\sD} d(z,gx_0)\]
		for all $g\in G$. Therefore $zu \approx_{2\sD} zv$ which implies $uo \approx_{2\sD} ov$. It follows that $o \approx_\sD m$, where $m$ is the midpoint of $[u,v]$. Finally, applying Lemma \ref{lemcirc2} gives
		$p \approx_{4\delta} o \approx_\sD m \approx_{3\delta} c$
		and we are done.

		\endproof

		\proofof{Proposition \ref{covercirc}}
		Let $\gamma'$ be a circumcentre for $G\alpha$. Combining Lemma \ref{lempishort} and Proposition \ref{hullshort}, we deduce 
		\[d_\Sigma(\pi(\alpha), \Hull(G\alpha)) \leq \sk_5',\]
		where $\sk_5'$ depends only on $\xi(\Sigma)$ and $\deg P$. Since $\pi(\alpha)$ is a $G$--invariant multicurve, the radius of its $G$--orbit is at most 1. Applying the previous lemma gives $d_\Sigma(\pi(\alpha), \gamma') \leq \sk_5' + 7\delta + 1$ and we are done.

		\halmos
		
		\subsection{Almost fixed point sets}
		
				Let $G$ be a finite group acting by isometries on a $\delta$--hyperbolic space $\X$. Given $\sR\geq 0$, let
		\[ \Fix_{\X}(G,\sR) := \left\{ x\in\X ~|~ \diam(Gx) \leq \sR \right\} \]
		be the set of \emph{$\sR$--almost fixed points} of $G$ in $\X$.

		\begin{Lem}\lab{lemfixrd}
		The set $\Fix_{\X}(G,2\delta)$ is non-empty. Moreover, for all $\sR\geq \delta$, 
		\[\Fix_{\X}(G,2\,\sR) \approx_{\sR + \delta} \Fix_{\X}(G,2\delta). \]
		
		\end{Lem}
		
		\proof
		Let $x$ be any point in $\X$ and let $c$ be a circumcentre for $Gx$. Since $Gx$ is $G$--invariant, all $G$--translates of $c$ are also circumcentres for $Gx$. It follows from Lemma \ref{lemcirc} that $c$ is contained in $\Fix_{\X}(G,2\delta)$.
		
		Assume further that $x$ is contained in $\Fix_{\X}(G,2\,\sR)$, where $\sR \geq \delta$.  By Lemma \ref{lemcirc2}, we have
		\[xc \leq \rad(Gx) \leq \frac{1}{2}\,\diam(Gx) + \delta \leq \sR + \delta \]
		and hence $\Fix_{\X}(G,2\,\sR) \subseteq_{\sR + \delta} \Fix_{\X}(G,2\delta)$. The reverse inclusion is immediate.
		\endproof
		
		Thus, to understand the geometry of the set of $2\sR$--fixed points of $\X$, for $\sR \geq \delta$, it suffices to study the behaviour of $\Fix_{\X}(G,2\delta)$. One can also check that $\Fix_{\X}(G,2\sR)$ is quasi-convex for $\sR\geq \delta$.
		
		\begin{Lem}\lab{lemfixpc}
		Let $\sR \geq \delta$. Given a point $x\in\X$, let $c$ be a circumcentre for its $G$--orbit. Let $p$ be a nearest point projection of $x$ to $\Fix_{\X}(G,2\sR)$. Then $c \approx_{\sk_6} p$, where $\sk_6 = 2\delta + 4\sR$. 
		
		\end{Lem}
		
		\proof
		For all $g\in G$, we have
		\[d(gx,p) \leq d(gx,gp) + d(gp,p) \leq d(x,p) + 2\sR \leq d(x,c) + 2\sR \leq \rad(Gx) + 2\sR.\]
		Applying Lemma \ref{lemcirc} completes the proof.
		\endproof
		
		We will demonstrate below that when $\sR < \delta$, Lemmas \ref{lemfixrd} and \ref{lemfixpc} need not hold; it is possible for $\Fix_{\X}(G,2\sR)$ to be empty or lie very deeply inside $\Fix_{\X}(G,2\delta)$.
		 
		 Recall that a point specified by ${(r,\theta,t) \in [0,\infty) \times \R \times \R}$ in the standard cylindrical co-ordinate system on $\R^3$ represents the point ${(r\cos\theta,r\sin\theta,t) \in\R^3}$ under the standard Cartesian co-ordinate system.
		 
		 \begin{Ex}[Rocketship]
		 A \emph{rocketship} of length $l>0$ with $n \geq 2$ fins, denoted $\rocket = \rocket(n,l)$, is the union of the following three subsets of $\R^3$ defined using cylindrical co-ordinates:
		  \begin{itemize}
				\item the \emph{nose} $\nose = \{(t,\theta,t) ~|~ 0\leq t \leq 1,\, \theta\in\R\}$, a right circular cone of height 1 and base radius 1;
				\item the \emph{shaft} $\shaft = \{(1,\theta,t) ~|~ 1 \leq t \leq l+1, \, \theta\in\R\}$, a right circular cylinder of height $l$ and base radius 1; and
				\item the \emph{fins} $\fins_n = \{(1,\frac{2k\pi}{n},t) ~|~ t \geq l+1,\, k\in\Z \}$, a disjoint union of $n$ closed rays.
			\end{itemize}
		We endow $\rocket$ with the path metric inherited from $\R^3$ equipped with the standard Euclidean metric. One can show that $\rocket$ is quasi-isometric to a tree and therefore $\delta$--hyperbolic for some $\delta>0$; this can be done by collapsing the radial component of the nose and shaft. Moreover, one can check that $\delta \geq \frac{\pi}{2}$ for $l$ sufficiently large.
		
		Observe that $G = \Z/n\Z$ acts isometrically on $\rocket$ by rotations about the {$t$--axis} through integral multiples of $\frac{2\pi}{n}$. Let $x$ be any point on $\fins_n$. Then a circumcentre $c$ for $Gx$ is a point of the form $(1,\frac{(4k+1)\pi}{2n},l+1)$, for some $k\in\Z$. For $\sR \geq 0$ sufficiently small, $\Fix_{\rocket}(G,2\sR)$ is contained entirely within the nose $\nose$. Therefore, $c$ must be a distance at least $l$ away from any nearest point projection of $x$ to $\Fix_{\rocket}(G,2\sR)$. Furthermore, $\Fix_{\rocket}(G,2\delta)$ contains both $\nose$ and $\shaft$, and so its Hausdorff distance from $\Fix_{\rocket}(G,2\sR)$ is at least $l$.
		 \end{Ex}
		 
		 Let us return our attention to the curve graph. When $G$ is the deck transformation group of a regular cover $P : \Sigma \rightarrow S$, the vertices in $\Fix_{\C(\Sigma)}(G,1)$ coincide exactly with those of $\Pi(\C(S)) \subseteq \C(\Sigma)$. 
		  Combining Theorem \ref{coverproj} and Proposition \ref{covercirc}, we deduce:
		  
		  \begin{Cor}
		  Any circumcentre for the $G$--orbit of a curve $\alpha\in\C(S)$ is a uniformly bounded distance away from any nearest point projection of $\alpha$ to $\Pi(\C(S))$.
		  \end{Cor}
		  
		   Therefore Lemma \ref{lemfixpc} still holds for $\Fix_{\C(\Sigma)}(G,1)$, albeit with weaker control over the constant $\sk_6$. As the example above shows, this cannot be proved using purely synthetic methods assuming only $\delta$--hyperbolicity of $\C(\Sigma)$. In conclusion: ``There are no rocketships in the curve complex.''

	\section{Singular Euclidean structures} \lab{chap-sing}
	
In this section we give a generalisation of Bowditch's singular Euclidean surfaces which are used to to estimate weighted intersection numbers in \cite{bhb-int}. We prove that such surfaces satisfy a quadratic isoperimetric inequality and so, by a theorem of Bowditch, must contain essential annuli of definite width, leading to a proof of Lemma \ref{boundshort}.

		\subsection{Construction of $S(\tvec\cdot\balpha)$}\lab{secbuild}

		Suppose $S = (S,\Omega)$ is a closed surface of genus $g$ with a set of $m$ marked points $\Omega$ so that $\xi(S) \geq 2$. Throughout this chapter, we shall fix an $n$--tuple $\balpha = (\alpha_1,\ldots,\alpha_n)$ of distinct multicurves in $\C(S)$ and a weight vector $\tvec = (t_1,\ldots,t_n)$. For simplicity, assume that $\balpha$ fills $S$ and that all entries of $\tvec$ are positive. We will deal with the appropriate modifications for the non-filling case in Section \ref{secnonfill}.

		Begin by realising the multicurves $\alpha_i$ on $S$ so that they intersect generally and pairwise minimally. We can achieve this, for example, by placing a complete hyperbolic metric on the complement of the marked points in $S$, taking geodesic representatives of the $\alpha_i$ and perturbing slightly if required. The union of the $\alpha_i$ gives a connected $4$--valent graph $\Upsilon$ on $S$. The closure of each component of $S - \Upsilon$ is a polygon with at most one marked point. The polygons together with $\Upsilon$ give $S$ the structure of a 2--dimensional cell complex. By taking the dual $2$--cell structure, we obtain a tiling of $S$ by rectangles which are in bijection with the self-intersection points of $\balpha$. We will insist that any marked point of $S$ coincides with a vertex of this tiling.

		Each rectangle $R$ corresponding to an intersection of $\alpha_i$ with $\alpha_j$ is isometrically identified with a Euclidean rectangle of side lengths $t_i$ and $t_j$ so that $\alpha_i$ is transverse to the two sides of length $t_i$. This gives a singular Euclidean metric on $S$. We may arrange for each $\alpha_i$ to be locally geodesic in this metric by requiring $\alpha_i \cap R$ to be a straight line connecting the midpoints of opposite sides of $R$, for every rectangle $R$ meeting $\alpha_i$. Thus, each component of $\alpha_i$ is the core curve of an annulus of width $t_i$ formed by taking the union of all rectangles $R$ it meets.

		The singular Euclidean surface defined above shall be denoted $S(\tvec\cdot\balpha)$. We remark that the metric depends on the realisation of $\balpha$ on $S$ up to isotopy, however, any such choice will work equally well for the purposes of proving the proposition below.

		We will allow ourselves to homotope a curve $\gamma\in\C(S)$ to meet marked points in order to speak of geodesic representatives on $S(\tvec\cdot\balpha)$. To be more precise, suppose $c'$ is a simple closed curve on $S$ representing $\gamma$. We say $c$ is a \emph{representative} of $\gamma$ if there is a homotopy $\bF : S^1 \times [0,1] \rightarrow S$ such that $\bF(\theta,0) = c'(\theta)$, $\bF(\theta,1) = c(\theta)$ and $\bF(S^1\times\{t\}) \subseteq S - \Omega$ for all $0\leq t < 1$.
		
			A geodesic representative $c$ of $\gamma$ on $S(\tvec\cdot\balpha)$ may not necessarily be an embedded copy of $S^1$. In these cases, there is a decomposition of the circle ${S^1 = \cup I_k}$ into a finite union of closed intervals with disjoint interiors so that $c : S^1 \rightarrow S(\tvec\cdot\balpha)$ embeds each $I_k$ as a straight line segment whose endpoints are singular points or marked points.
			
			We will use $l(\gamma)$ to denote the length of a geodesic representative of $\gamma$ on $S(\tvec\cdot\balpha)$ with respect to the singular Euclidean metric.

				\begin{Prop}\lab{sing}
		The singular Euclidean surface $S(\tvec\cdot\balpha)$ has the following properties.
			\begin{enumerate}
			 \item $S(\tvec\cdot\balpha)$ has area $\norma^2 = \sum_{j<k} t_jt_ki(\alpha_j,\alpha_k)$.
			 \item For all curves $\gamma\in\C(S)$, we have
				\[l(\gamma) \leq i(\tvec\cdot\balpha,\gamma) \leq \sqrt{2}l(\gamma).\]
			 \item There exists an essential annulus on $S(\tvec\cdot\balpha)$ whose width is at least $\sW_0 \norma$, where $\sW_0 > 0$ is a constant depending only on $\xi(S)$.
			\end{enumerate}
				\end{Prop}

		The first claim is immediate from the construction. Before proving the second and third claims, we will outline some consequences of this proposition which shall later be used to prove  Lemma \ref{boundshort}. It is worth mentioning that the third claim holds for a larger class of metrics satisfying a suitable isoperimetric inequality. We also remark that the metric on $S(\tvec\cdot\balpha)$ can be approximated by a non-singular Riemannian metric. We will, however, choose to work with singular Euclidean metrics to simplify the exposition.

\subsection{Short curves and wide annuli}\lab{secannuli}
		
		We now state some facts arising from the interplay between weighted intersection numbers, lengths of curves and widths of annuli on $S(\tvec\cdot\balpha)$. Most of the statements and proofs are covered in \cite{bhb-int} but we will include them for completeness.
		
		Let $A$ be a closed Riemannian annulus. Define $\width(A)$ to be the length of a shortest arc joining its two boundary components and $\length(A)$ to be the length of a shortest core curve of $A$. The following is a consequence of the Besicovitch Lemma \cite{Besicovitch} (see Lemma 4.5$\frac{1}{2}$ in \cite{Gromov-metric} for a proof).
		
		\begin{Lem}[Annulus inequality]\lab{lemannulus}
		Let $A$ be an annulus. Then
		\[ \width(A)\times\length(A) \leq \area(A) \rlap{\qquad \qquad \qquad \qquad \qquad \halmos}\]
		\end{Lem}

		\noindent If $\gamma$ is the core curve of an annulus $A$ on $S(\tvec\cdot\balpha)$ then
		\[ \width(A)\times l(\gamma) \leq \width(A)\times\length(A) \leq \area(A) \leq \area(S(\tvec\cdot\balpha)) = \norma^2\]
		where we have applied the annulus inequality for the second comparison.

		\begin{Lem}\lab{annul}
		Let $A$ be an annulus on $S(\tvec\cdot\balpha)$ with core curve $\gamma$. Then for all $\beta\in\C(S)$, we have
		\[ \width(A)\times i(\gamma,\beta) \leq l(\beta).\]
		\end{Lem}
		
		\proof
		Let $b ~:~ S^1 \rightarrow S(\tvec\cdot\balpha)$ be an map which realises $\beta$ as a geodesic on $S(\tvec\cdot\balpha)$. We may pull back the metric on $S(\tvec\cdot\balpha)$ via $b$ to give $S^1$ the structure of a circle of length $l(\beta)$. The result follows by observing that $b^{-1}(A)\subseteq S^1$ contains at least $i(\gamma,\beta)$ disjoint arcs, each having length at least $\width(A)$.
		\endproof

		By combining the above inequalities with Proposition \ref{sing}, we obtain the following bounds on intersection number with the core curve of an annulus: 
		
		\begin{Cor}\lab{annul2}
		Let $A$ be an essential annulus on $S(\tvec\cdot\balpha)$ with core curve $\gamma$. Then
		\[i(\tvec\cdot\balpha,\gamma) \leq \frac{\sqrt{2}\norma^2}{\width(A)} \textrm{~ and~ } i(\gamma,\beta) \leq \frac{i(\tvec\cdot\balpha,\beta)}{\width(A)}\]
		for all $\beta\in\C(S)$. \halmos

		\end{Cor}

		\subsection{A grid structure on $S(\tvec\cdot\balpha)$}
		
			In this section, we describe a grid structure on $S(\tvec\cdot\balpha)$ and give a proof of the second claim of Proposition \ref{sing}.
		
		\begin{Dfn}[Quarter-translation surface]
		A \emph{quarter-translation surface} is a topological surface $S$ with a finite set of singularities $\varsigma$ together with an atlas of charts from $S - \varsigma$ to $\R^2$ whose transition maps are translations of $\R^2$ possibly composed with rotations through integral multiples of $\frac{\pi}{2}$.
		
		The singular points have cone angles which are integral multiples of $\frac{\pi}{2}$ and at least $\pi$.
		\end{Dfn}

		Suppose $S$ is a quarter-translation surface. We may pull back the standard Euclidean metric on $\R^2$ to give a singular Euclidean metric on $S$. Geodesics which do not meet any singular points or marked points with respect to this metric can only self-intersect orthogonally. We will also define an $L^1$ metric on $S$ by pulling back the metric given infinitesimally by $|dx| + |dy|$ on $\R^2$. The following is immediate:

		\begin{Lem}
		Let $l^2(\eta)$ and $l^1(\eta)$ denote, respectively, the Euclidean and $L^1$ lengths of a path $\eta$ on $S$. Then $l^2(\eta) \leq l^1(\eta) \leq \sqrt{2} l^2(\eta)$. \halmos
		\end{Lem}	
	
		Whenever we deal with quarter-translation surfaces, we will assume that we are working with the singular Euclidean metric unless otherwise specified.		

		The transition maps between the charts for $S$ preserve a pair of orthogonal directions on $\R^2$ which we may take to be the standard horizontal and vertical directions. By pulling these back via the coordinate charts, we can equip $S$ with a preferred (unordered) pair of orthogonal directions defined on the complement of the singular points. We shall refer to these as the \emph{grid directions}. Geodesics which run parallel to a grid direction will be called \emph{grid arcs}. Every non-singular point on $S$ has an open rectangular neighbourhood, with sides parallel to the grid directions, on which the grid leaves restrict to give a pair of transverse foliations. Such a rectangle will be called an \emph{open grid rectangle}.

		It is straightforward to check that $S(\tvec\cdot\balpha)$ is a quarter-translation surface. We will assume that the grid directions on $S(\tvec\cdot\balpha)$ run parallel to the sides of the rectangles used in its construction.

		\begin{Lem}
		Given a curve $\gamma\in\C(S)$, let $c$ be any of its geodesic representatives on $S(\tvec\cdot\balpha)$ with respect to the Euclidean metric. Then
		\[l^1(c) = i(\tvec\cdot\balpha,\gamma). \]
		\end{Lem}
		
		\proof
		
		If $c$ is an embedded simple closed loop then we can isotope it to another geodesic representative which meets at least one singularity. Thus we can assume that $S^1$ decomposes as a finite union of intervals $\cup I_k$ with disjoint interiors such that $c~:~ S^1 \rightarrow S(\tvec\cdot\balpha)$ embeds each $I_k$ as a straight line segment with singularities or marked points at its endpoints.
		
		We can homotope $c$ to a closed path $c' ~:~ S^1 \rightarrow S(\tvec\cdot\balpha)$ so that each $c'(I_k)$ is an edge--path in the {1--skeleton} of $S(\tvec\cdot\balpha)$ with the same endpoints as $c(I_k)$. The homotopy can be performed in a way which preserves the $l^1$--length of the path and without creating new intersection points with any of the $\alpha_i$. One can check that $c$ intersects each $\alpha_i$ minimally and thus the same is also true of $c'$. Finally, we deduce
		\[l^1(c'(S^1)) = i(\tvec\cdot\balpha,\gamma)\]
		by observing that every edge in the 1--skeleton of $S(\tvec\cdot\balpha)$ transverse to $\alpha_i$ has length $t_i$.
		\endproof

		The second claim of Proposition \ref{sing} follows from the previous two lemmas.

		\subsection{An isoperimetric inequality}\lab{secisop}

		Let $S = (S,\Omega)$ be a closed singular Riemannian surface $S$ with a finite set of marked points $\Omega$. Let $\Delta$ be a closed disc and suppose $\iota: \Delta \rightarrow S$ is a piecewise smooth immersion which restricts to an embedding on its interior. Let $D$ denote the image $\iota(\inte(\Delta))$.

		\begin{Dfn}
		An open disc $D$ arising in the above manner is called a \emph{trivial region} on $S$ if it contains at most one marked point.
		\end{Dfn}

		Bowditch defines trivial regions as open discs on $S$ containing at most one marked point without any conditions concerning piecewise smooth embeddings. Nevertheless, his proof of the following proposition still holds with our definition:
		
		\begin{Prop}[\cite{bhb-int}]\lab{propisop}
		Suppose $f ~:~ [0,\infty) \rightarrow [0,\infty)$ is a homeomorphism. Let $\rho$ be a singular Riemannian metric on an orientable closed surface $S$ with unit area. Let $\Omega$ be a finite set of marked points on $S$. We will assume $|\Omega|\geq 5$ whenever $S$ is a 2--sphere. If $\area(D) \leq f(\length(\partial D))$ for any trivial region $D$ then there is an essential annulus $A \subseteq S- \Omega$ such that $\width(A) \geq \sW_0$, where $\sW_0 > 0$ depends only on $\xi(S)$ and $f$.
		\end{Prop}

		A little care is required to clarify what $\length(\partial D)$ means, especially when $\partial D$ is not an embedded copy of a circle. In general, the boundary $\partial D$ is an embedded Eulerian graph on $S$ whose edges are piecewise smooth arcs. We define $\length(\partial D)$ to be the sum of the lengths of the arcs with respect to the metric on $S$.

		This section will be devoted to proving the following lemma which, together with the above proposition, implies the third claim of Proposition \ref{sing}.

		\begin{Lem}\lab{lemisop}
		Suppose $D$ is a trivial region on $S(\tvec\cdot\balpha)$. Then
		\[\area(D) \leq 4\,\length(\partial D)^2.\]
		\end{Lem}

		Before launching into the details of the proof, we give a brief outline of our argument. First, we reduce the problem to that of studying embedded closed discs on $S(\tvec\cdot\balpha)$ whose boundary is a finite union of grid arcs. We then show that such a disc $D$ can be given a tiling by grid rectangles. This tiling is dual to a collection of arcs on $D$, where each arc is parallel to a component of some $\alpha_i \cap D$. We shall call the union of all rectangles meeting a given arc a \emph{band}. The key step is to observe that any two arcs in the collection intersect at most twice. Thus, the intersection of two distinct bands is the union of at most two rectangles arising from the tiling. Conversely, any rectangle from the tiling is contained in the intersection of two such bands. This then allows us to bound the area of the rectangles in terms of the length of $\partial D$.
		
		\subsubsection{Technical adjustments}

		Let us first make a couple of observations to simplify the problem.

		\begin{Lem}\lab{lemgrid}
		Any trivial region $D$ on $S(\tvec\cdot\balpha)$ can be perturbed to a trivial region $D'$ whose boundary is a finite union of grid leaves. Moreover, $D'$ can be chosen so that $\area(D') \geq \area(D)$ and $\length(\partial D') \leq \sqrt{2} ~\length (\partial D)$. \halmos
		\end{Lem}

		We will henceforth assume that the boundary of any trivial region on $S(\tvec\cdot\balpha)$ is a finite union of grid leaves.

		Let $\iota : \Delta \rightarrow S(\tvec\cdot\balpha)$ be a piecewise smooth immersion whose restriction to $\inte (\Delta)$ is an embedding with image $D$. Observe that $\iota: \partial \Delta = S^1 \rightarrow \partial D$ is an immersion of a circle which runs over each edge of $\partial D$ at most twice. We will metrise $\Delta$ by pulling back the metric on $S(\tvec\cdot\balpha)$ via $\iota$. 

		\begin{Lem}\lab{lemdisc}
		Suppose $D$ and $\Delta$ are as given above. Then $\area(\Delta) = \area(D)$ and $\length(\partial D) \leq \length(\partial \Delta) \leq 2~\length(\partial D)$. \halmos
		\end{Lem}
		
		\subsubsection{Tiling $\Delta$ by rectangles}

		The disc $\Delta$ inherits grid directions from $S(\tvec\cdot\balpha)$ via $\iota$ away from the preimage of the singular points. The boundary decomposes as a finite union $\partial \Delta = \cup I_k$ of closed grid arcs with disjoint interiors. We may assume that this decomposition is minimal, that is, it cannot be obtained from any other such decomposition by subdividing arcs. An endpoint of any grid arc $I_k$ will be called a \emph{corner point} of $\partial \Delta$. A corner point which does not coincide with a singularity or a marked point must be an orthogonal intersection point of two grid arcs.

		It is worth noting that $\partial \Delta$ must contain at least two corner points and at least three if $D$ contains no marked points. To see this, recall that the grid leaves on $S(\tvec\cdot\balpha)$ are parallel to some $\alpha_i$. Any of the forbidden cases will imply that some $\alpha_i$ is trivial, peripheral, self-intersects or does not intersect some $\alpha_j$ minimally.

		Let us refer to marked points, corner points and singularities collectively as \emph{bad points}. Let $Z\subset\Delta$ be the union of $\partial\Delta$ with all grid arcs in $\Delta$ which have a bad point for at least one of their endpoints. Since there are finitely many bad points in $\Delta$, it follows that $Z$ is a finite embedded graph on $\Delta$. A vertex $v\in\inte\Delta\cap Z$ has valence $k$ if and only if the cone angle at $v$ is $\frac{k\pi}{2}$. If $v$ is a vertex which lies on $\partial\Delta$ then it has valence $k+1$ if and only if the cone angle at $v$ inside $\Delta$ is $\frac{k\pi}{2}$. It follows that every vertex $v$ of $Z$ has valence at least 2, and at least 3 if $v$ is not a marked point.

		\begin{Lem}\lab{lemtilegrid}
		There exists a tiling of $\Delta$ by finitely many grid rectangles with $Z$ as its 1--skeleton.
		\end{Lem}
		
		\proof
		First note that there are finitely many connected components of $\Delta- Z$ since $Z$ is a finite graph. Let $R$ be such a component and let $\bar R$ be its completion  with respect to its induced path metric. Observe that $\bar R$ is a closed planar region admitting a Euclidean metric with piecewise geodesic boundary, where the interior angle between adjacent edges of $\partial \bar R$ is $\frac{\pi}{2}$. By the Gauss--Bonnet formula, the sum of its interior angles must equal $2\pi \chi(R)$. Since the frontier of $R$ in $\Delta$ meets at least one vertex of $Z$, the angle sum must be strictly positive. As $R$ is planar, it follows that $\chi(R) = 1$ and therefore $\bar R$ is a Euclidean rectangle. Also note that $Z$ is connected, for otherwise there would exist some component of $\Delta- Z$ with disconnected frontier.
		
		The inclusion $R \hookrightarrow \Delta$ can be extended continuously to a map $\bar R \rightarrow \Delta$, sending each edge of $\partial \bar R$ isometrically to an edge of $Z$ meeting the frontier of $R$. Thus $R$ is a grid rectangle since the edges of $Z$, by construction, are parallel to the grid directions. Finally, the closures of distinct rectangles $R$ and $R'$ can only intersect in a union of vertices and edges of $Z$.	
	
		\endproof

	\subsubsection{Controlling the area}

	Let $\A$ be the set of maximal grid arcs in $\Delta$ which intersect $Z$ only at midpoints of edges of $Z$. This is a collection of arcs dual to the tiling of $\Delta$ by rectangles described in Lemma \ref{lemtilegrid}. (There cannot exist any closed curves in $\Delta$ dual to the tiling as this would imply that some $\alpha_i$ is trivial or peripheral.) Given an arc $a\in\A$, let $B = B(a)$ be the union of all rectangles in the tiling which meet $a$. We will call $B$ a \emph{band} and $a$ a \emph{core arc} of $B$. Define $\width(B)$ to be the length of any edge of $Z$ crossed by $a$. Note that the set of bands in $\Delta$ is in bijection with $\A$.
	
	\begin{Lem}\lab{lemrect}
	The intersection of two distinct bands $B$ and $B'$ is the union of at most two rectangles whose side lengths are $\width(B)$ by $\width(B')$. Conversely, each rectangle in the tiling lies in the intersection of a unique pair of distinct bands.
	\end{Lem}
	
	\proof

	Let $a$ and $a'$ be core arcs of $B$ and $B'$ respectively. If $a$ and $a'$ intersect at least 3 times then they must bound a bigon in $\Delta$ containing no marked points. Now, $a$ and $a'$ can both be properly isotoped in $\Delta$ to components of ${\iota^{-1}(\alpha_i \cap D)}$ and ${\iota^{-1}(\alpha_i \cap D)}$ for some $\alpha_i$ and $\alpha_j$ respectively. Moreover, the isotopies can be performed without passing through any singular points or marked points. This procedure cannot destroy any bigons since any right-angled bigon on $\Delta$ must contain at least one singularity. It follows that $\alpha_i$ and $\alpha_j$ also bound a bigon in $D$, contradicting minimality.

	For the converse, simply take the bands corresponding to the unique pair of arcs which have an intersection point inside the given rectangle.
	\endproof
	
	We will refer to an edge of $Z$ lying in $\partial \Delta$ simply as an \emph{edge} of $\partial \Delta$.
	
	\begin{Lem}\lab{lembandwidth}
	Under the above hypotheses, we have
	\[\length(\partial \Delta) = 2\sum_B \width(B)\]
	where the sum is taken over all bands $B$ in $\Delta$.
	\end{Lem}
	
	\proof
	Let $B$ be a band in $\Delta$. Observe that $B \cap \partial\Delta$ consists of exactly two edges of $\partial\Delta$ whose length is equal to $\width(B)$. Conversely, each edge of $\partial\Delta$ lies in exactly one band.
	\endproof

		\begin{Lem}
		Let $\Delta$ be as above. Then
		\[\area(\Delta) \leq \frac{1}{2}~\length(\partial \Delta)^2. \]
		\end{Lem}

		\proof
		By Lemma \ref{lemrect}, $\Delta$ is a union of rectangles, each of which lies in the intersection of a pair of distinct bands. Thus
		\[\area(\Delta) =  \area\left(\bigcup_{B \neq B'} B \cap B'\right) = \sum_{B \neq B'} \area(B \cap B').\]
		Since the intersection of two distinct bands is the union of at most two rectangles whose side lengths are equal to the widths of the bands, we have
		\[\area(B \cap B') \leq 2~\width(B)\times\width(B')\]
		and hence
		\[\area(\Delta) \leq 2~\sum_{B \neq B'} \width(B)\times\width(B') \leq 2 ~\left(\sum_{B} \width(B)\right)^2.\]
		Finally, applying Lemma \ref{lembandwidth} completes the proof.
		\endproof

		Combining this result with Lemmas \ref{lemgrid} and \ref{lemdisc} completes the proof of Lemma \ref{lemisop}. \halmos
		
		\subsection{Non-filling curves}\lab{secnonfill}

			We now generalise the construction of $S(\tvec\cdot\balpha)$ to encompass non-filling curves. Assume $\balpha = (\alpha_1,\ldots,\alpha_n)$ is an $n$--tuple of distinct multicurves and $\tvec=(t_1,\ldots,t_n) \neq \mathbf{0}$ is a weight vector satisfying $\norma \neq 0$. Realise $\balpha$ minimally on $S$ to form a {4--valent} graph $\Upsilon$ on $S$. 

	Let $\Sigma\subseteq S$ be the (possibly disconnected) subsurface filled by $\Upsilon$. This can be obtained by taking a closed regular neighbourhood of $\Upsilon$ on $S$ and then attaching all complementary regions which are discs with at most one marked point.		
			If $\balpha$ fills $S$ then $\Sigma = S$. In general, $\Sigma$ will be a disjoint union of surfaces $\Sigma_1 \cup \ldots \cup \Sigma_s$. Observe that $s \leq \xi(S)$ since we can find a multicurve on $S$ so that exactly one component is contained in each $\Sigma_k$ (by choosing a suitable subset of all curves appearing in $\balpha$, for example). Some of these components may be annuli -- this occurs precisely when a multicurve $\alpha_i$ has a component disjoint from all other $\alpha_j$. All other components will have genus at least one, or are spheres where the sum of the number of marked points and boundary components is at least four.

			We now define a 2--dimensional complex $S(\tvec\cdot\balpha)$ as a quotient of $S$. Suppose $\Sigma_k$ is an annular component of $\Sigma$ whose core curve is a component of $\alpha_i$. We identify $\Sigma_k$ with $S^1 \times [0,t_i]$ and then collapse the first co-ordinate to give a closed interval $I_k$ of length $t_i$. Next, we collapse every complementary component of $\Sigma$ in $S$ to a marked point. These marked points will be called \emph{essential}. We then apply the construction given in Section \ref{secbuild} to the image of each non-annular component of $\Sigma$ in the quotient space. The resulting space is a finite collection of singular Euclidean surfaces and closed intervals identified along appropriate essential marked points. Note that this construction agrees with the one given in Section \ref{secbuild} for the case of filling curves.
			
			Let $c$ be a representative of a curve $\gamma\in\C(S)$ on $S$. Its image $\bar c$ on $S(\tvec\cdot\balpha)$ will be a closed curve or a union of paths connecting essential marked points. Define $l(\gamma)$ to be the minimal length of $\bar c$ over all representatives $c$ of $\gamma$.

			\begin{Prop}\lab{sing2}
			Suppose $\balpha$ and $\tvec$ satisfy $\norma > 0$. Then the first two claims of Proposition \ref{sing} hold for $S(\tvec\cdot\balpha)$. \halmos
			\end{Prop}
			
			The proof of the above proceeds in the same manner as for the case of filling curves. It remains to prove an analogue of the third claim.

		We will refer to the image of a component $\Sigma_k$ of $\Sigma$ in $S(\tvec\cdot\balpha)$ as a \emph{component} of $S(\tvec\cdot\balpha)$. Let $Y$ be a component of $S(\tvec\cdot\balpha)$ with maximal area. Since $\Sigma$ has at most $\xi(S)$ components, we have $\area(Y) \geq \frac{\norma^2}{\xi(S)}$. Note that $Y$ cannot be an interval since $\norma > 0$.

		We may argue as in Section \ref{secisop} to show the following.
		
		\begin{Lem}\lab{lemnonfill}
		Suppose $Y$ has genus at least one, or is a sphere with at least five marked points. Then there is an essential annulus on $Y$ whose width is at least $\frac{\sW_0\norma}{\sqrt{\xi(S)}}$, where $\sW_0$ is a constant depending only on $\xi(Y) \leq \xi(S)$. \halmos
		\end{Lem}

		Let us assume $Y$ is a sphere with 4 marked points for the rest of this section. We will not prove the existence of annuli of definite width on $Y$. Instead, we show that it suffices to find wide annuli on a torus $\TT^2 = \R^2 / \Z^2$ which double branch covers $Y$ for the purposes of proving Lemma \ref{boundshort}.

		Define a hyperelliptic involution $h : \TT^2 \rightarrow \TT^2$ by setting $h(x,y) = (-x,-y)$ modulo $\Z^2$. The quotient map induced by the action of $\ang{h}$ is a double cover ${P: \TT^2 \rightarrow S^2}$ branched over four points. The branch points correspond to the fixed orbits $(0,0)$, $(0,\frac{1}{2})$,  $(\frac{1}{2},0)$ and $(\frac{1}{2},\frac{1}{2})$ of $h$.
		
		We will identify $Y$ with the quotient space $S^2$ so that the marked points coincide exactly with the branch points. We then metrise $\TT^2$ by lifting the singular Euclidean metric on $Y$ via $P$.
		This metric can also be obtained by taking the preimages of the $\alpha_i$ contained in $Y$ to $\TT^2$ and then applying the construction as described in Section \ref{secbuild}. It follows from the work in Section \ref{secisop} that $\TT^2$ enjoys the isoperimetric inequality stated in Lemma \ref{lemisop}. Invoking Proposition \ref{propisop} gives the following:
		
		\begin{Lem}\lab{lemtoruswide}
		There exists an essential annulus on $\TT^2$ of width at least $\frac{\sW\norma}{\sqrt{\xi(S)}}$ for some universal constant $\sW>0$. \halmos
		\end{Lem}
		
		\begin{Rem}
		By following Bowditch's proof of Proposition \ref{propisop} in \cite{bhb-int} for the the case of the torus, one can show that setting $\sW = \frac{1}{3\sqrt{2}}$ will satisfy the statement of the above lemma.
		\end{Rem}

		Observe that $h(\tilde \gamma)$ is homotopic to $\tilde \gamma$ for any simple closed curve $\tilde \gamma$ on $\TT^2$. Thus, the image of any simple closed on $\TT^2$ under $P$ is homotopic to a simple closed curve on $Y$.
		
		\begin{Lem}\lab{lemtoruswide2}
		Let $A$ be an essential annulus on $\TT^2$ with core curve $\tilde \gamma$. Let $\gamma\in\C(S)$ be the image of $\tilde \gamma$ on $Y$ under $P$. Then
		
		\[i(\gamma,\beta) \leq \frac{2\, i(\tvec\cdot\balpha,\beta)}{\width(A)}\]
		for all $\beta\in\C(S)$.
		\end{Lem}
		
		\proof

		First, recall that $\beta \cap Y$ is either a simple closed curve or a union of paths connecting marked points of $Y$. The preimage $P^{-1}(\beta)$ is a finite union of (not necessarily disjoint) essential curves on $\TT^2$. By perturbing $\gamma$ to an embedded curve which misses the marked points of $Y$, we see that each point of $\gamma \cap \beta$ lifts to exactly two points on $\TT^2$ under $P$, and so
		\[i(\gamma,\beta) = \frac{i(P^{-1}(\gamma),P^{-1}(\beta))}{2} \leq i(\tilde \gamma,P^{-1}(\beta)).\]
		Applying Lemma \ref{annul} to each curve in $P^{-1}(\beta)$ gives
		\[\width(A) \times i(\tilde\gamma, P^{-1}(\beta)) \leq l(P^{-1}(\beta)). \]
		Next, we have
		\[l(P^{-1}(\beta)) = 2\, l(\beta \cap Y) \leq 2\, l(\beta) \leq 2\, i(\tvec\cdot\balpha,\beta)\]
		where we have applied Proposition \ref{sing2} for the final comparison.
		Finally, combining the above inequalities gives the desired result.

		\endproof

		\subsection{Proof of Lemma \ref{boundshort}} \lab{secproofshort}
		
		We finally show that
		\[\short(\tvec\cdot\balpha, \sL) = \left\{\gamma\in\C(S) ~|~ i(\tvec\cdot\balpha,\gamma) \leq \sL \norma \right\}\]
		has uniformly bounded diameter in $\C(S)$.

		If $\norma = 0$ then $\short(\tvec\cdot\balpha, \sL)$ contains the $\alpha_i$ and is contained in the {1--neighbourhood} of $\balpha$ in $\C(S)$.
		
		 Next, suppose $\norma > 0$ and $\balpha$ fills $S$. By the third claim of Proposition \ref{sing}, there exists an essential annulus $A$ on $S(\tvec\cdot\balpha)$ whose width is at least $\sW_0 \norma$, where $\sW_0$ depends only on $\xi(S)$. Set $\sL_0 = \frac{\sqrt{2}}{\sW_0}$. Let $\gamma$ be the core curve of $A$. By Corollary \ref{annul2}, we have 
		\[i(\tvec\cdot\balpha, \gamma) \leq \frac{\sqrt{2}\norma^2}{\width(A)} \leq \frac{\sqrt{2}}{\sW_0} \norma\]
		and hence $\short\left(\tvec\cdot\balpha,\sL \right) \neq \emptyset$ for all $\sL\geq \sL_0$.		
		Furthermore, if $\beta\in\short\left(\tvec\cdot\balpha,\sL \right)$ is another curve then, by Corollary \ref{annul2}, we have
		\[i(\gamma,\beta) \leq \frac{i(\tvec\cdot\balpha,\beta)}{\width(A)} \leq \frac{\sL\norma}{\sW_0\norma} = \frac{\sL}{\sW_0}. \]
		Applying Lemma \ref{lemhempel} and the triangle inequality gives
		\[\diam(\short(\tvec\cdot\balpha, \sL)) \leq 2\left[2\log_2 \left(\frac{\sL}{\sW_0}\right) + 2\right] = 4\log_2 \sL + \sk_0\]
		where $\sk_0$ is a constant depending only on $\xi(S)$.
		
		For the case where $\norma > 0$ but $\balpha$ does not fill $S$, it is immediate that ${\short(\tvec\cdot\balpha, \sL)}$ is non-empty. To bound the diameter, invoke Lemmas \ref{lemnonfill}, \ref{lemtoruswide} and \ref{lemtoruswide2} then argue as above.\halmos

	\bibliography{mybib}		
	\bibliographystyle{amsalpha}

\end{document}